\documentclass[11pt]{amsart}
\usepackage[T1]{fontenc}
\usepackage{lmodern,amsmath,amssymb,mathtools,microtype,mathrsfs}
\usepackage[margin=1in]{geometry}
\usepackage{xcolor}
\usepackage[colorlinks=true,linkcolor=blue,citecolor=green!45!black,urlcolor=blue!45!black]{hyperref}
\hypersetup{pdftitle={Accessibility of collapsed Anosov flows},pdfauthor={Ziqiang Feng and Ra\'ul Ures},pdfsubject={Accessibility in dimension three without recurrence assumptions}}
\newtheorem{theorem}{Theorem}[section]
\newtheorem{proposition}[theorem]{Proposition}
\newtheorem{lemma}[theorem]{Lemma}
\newtheorem{corollary}[theorem]{Corollary}
\theoremstyle{definition}

\theoremstyle{remark}

\ifcsname conjecture\endcsname\else
  \begingroup
  \theoremstyle{plain}
  \newtheorem{conjecture}{Conjecture}
  \endgroup
\fi

\newcommand{\R}{\mathbb R}
\newcommand{\Z}{\mathbb Z}
\newcommand{\cF}{\mathcal F}
\newcommand{\cW}{\mathcal W}
\newcommand{\cO}{\mathcal O}
\newcommand{\wt}[1]{\widetilde{#1}}
\newcommand{\Inv}{\operatorname{Inv}}
\newcommand{\id}{\operatorname{id}}
\newcommand{\Int}{\operatorname{int}}

\newcommand{\eps}{\varepsilon}

\title[Accessibility of collapsed Anosov flows]{Accessibility of collapsed Anosov flows}
\author{Ziqiang Feng}
\address{Beijing International Center for Mathematical Research, Peking University, Beijing 100871, China}
\email{zqfeng@pku.edu.cn}
\author{Ra\'ul Ures}
\address{Department of Mathematics and Shenzhen International Center for Mathematics, Southern University of Science and Technology, Shenzhen 518000, China}
\email{ures@sustech.edu.cn}
\date{\today}
\begin{document}

\begin{abstract}
We prove that every $C^1$ collapsed Anosov flow on a closed three-manifold with non-virtually-solvable fundamental group is accessible. No transitivity or nonwandering assumption is required, and the collapsing map need not satisfy the strong-collapse conditions. 
We also prove that every three-dimensional collapsed Anosov flow for which no positive iterate is a discretized Anosov flow is accessible.
\end{abstract}
\maketitle

\tableofcontents

\section{Introduction}

\subsection{Motivation and historical remarks}

The study of accessibility grew out of the effort to extend the Hopf argument beyond uniform hyperbolicity. For a partially hyperbolic diffeomorphism $f$ of a closed manifold $M$, the invariant splitting
\[
TM=E^s\oplus E^c\oplus E^u
\]
provides uniformly contracting and expanding directions, separated by a center direction with intermediate behavior. The strong bundles integrate uniquely to foliations $\cW^s$ and $\cW^u$. The map $f$ is \emph{accessible} if any two points of $M$ can be joined by a finite concatenation of arcs in these foliations. Accessibility thus describes how the two hyperbolic directions connect the manifold, even though they need not span its tangent bundle.

Brin and Pesin~\cite{BP74} recognized its role in the ergodic theory of partially hyperbolic systems. The stable ergodicity of the time-one map of a geodesic flow, established by Grayson, Pugh and Shub~\cite{GPS94} for surfaces of constant negative curvature, showed that this mechanism could persist under perturbations without smoothness of the invariant distributions. Subsequent work of Wilkinson~\cite{Wil98} and Burns, Pugh and Wilkinson~\cite{BPW00} extended the theory to variable curvature and to volume-preserving Anosov flows in dimension three.

These developments led Pugh and Shub~\cite{PS96,PS00} to conjecture that stable ergodicity is dense among volume-preserving partially hyperbolic diffeomorphisms. Their program separates the abundance of accessibility from its ergodic consequences. The accessibility part of their conjecture can be stated as follows.

\begin{conjecture}\label{conj:ps-accessibility}(Pugh-Shub)
For every closed connected manifold $M$ and every $r\ge2$, accessibility holds on a $C^1$-open and $C^r$-dense subset of the space of $C^r$ partially hyperbolic diffeomorphisms of $M$. 
\end{conjecture}

The other part of the Pugh--Shub program asks whether essential accessibility implies ergodicity in the conservative setting, where essential accessibility means that every measurable union of accessibility classes has zero or full volume. The geometric question makes sense independently of an invariant volume and remains equally relevant in the nonconservative setting.

Dolgopyat and Wilkinson~\cite{DW03} proved the $C^1$-density of stable accessibility in arbitrary center dimension, in both the conservative and nonconservative categories. When the center is one-dimensional, Didier~\cite{Did03} established the stability of accessibility, Hertz, Hertz and Ures~\cite{HHU08} proved its $C^r$-density in the conservative category, and Burns, Hertz, Hertz, Talitskaya and Ures~\cite{BRHTU08} obtained the corresponding nonconservative result for $r\ge2$. The recent work of Ures and Wang~\cite{UW26} strengthens this density to be $C^{\omega}$ in the real analytic setting. Thus the accessibility part of the Pugh--Shub program is settled for one-dimensional center. On the measure-theoretic side, Burns and Wilkinson~\cite{BW10} proved that essential accessibility implies ergodicity for $C^2$ volume-preserving systems under center bunching, without assuming dynamical coherence. This bunching condition is automatic in center dimension one, a case also treated by Hertz, Hertz and Ures~\cite{HHU08}. In arbitrary center dimension, following the approach of~\cite{HHTU11}, Avila, Crovisier and Wilkinson~\cite{ACW21} subsequently proved the $C^1$-density of stable ergodicity among volume-preserving $C^r$ partially hyperbolic diffeomorphisms, for $r>1$. See also~\cite{Wil10,CHHU18} for broader accounts of these developments.

The abundance of accessibility and stable ergodicity in dimension three leaves open the problem of determining which individual systems can fail to have these properties. Hertz, Hertz and Ures~\cite{HHU08b} proposed that the topology of the ambient manifold should sharply restrict the exceptional systems. The Hertz-Hertz-Ures Ergodicity Conjecture predicts that a nonergodic conservative partially hyperbolic diffeomorphism of a closed orientable three-manifold can occur only on a short list of torus bundles. The stronger formulation predicts an embedded torus tangent to $E^s\oplus E^u$ whenever ergodicity fails. Such a torus is an obstruction to accessibility, and the classification of Anosov tori~\cite{HHU11} explains the corresponding restriction on the ambient manifold. The Pugh--Shub program concerns an open and dense family of systems, whereas the Hertz--Hertz--Ures conjecture seeks to identify the geometric obstructions in the exceptional family. In particular, on manifolds with non-virtually-solvable fundamental group, it predicts ergodicity for every $C^2$ volume-preserving partially hyperbolic diffeomorphism.

The geometric nature of accessibility suggests a nonconservative strengthening of this picture. 
In the setting relevant here, this leads to the following formulation.

\begin{conjecture}\label{conj:hhu-accessibility}(Hertz-Hertz-Ures)
Let $f:M\to M$ be a $C^1$ partially hyperbolic diffeomorphism of a closed connected three-manifold. If $\pi_1(M)$ is not virtually solvable, then $f$ is accessible.
\end{conjecture}

We refer to this statement as the Hertz--Hertz--Ures Accessibility Conjecture. It asks for accessibility of each individual system, without an invariant volume or a recurrence hypothesis. An early precedent for accessibility forced by topology without conservativity or recurrence is the theorem of Hammerlindl and Potrie~\cite[Section 6.3]{HP14} for non-toral three-dimensional nilmanifolds. Hammerlindl and Shi~\cite{HS21} proved that every $C^{1+\alpha}$ partially hyperbolic diffeomorphism of $\mathbb T^3$ homotopic to an Anosov automorphism is accessible unless it is itself Anosov. In the non-virtually-solvable setting, Fenley and Potrie~\cite{FP22} proved the Conjecture above for discretized Anosov flows. Our recent work~\cite{FU26b} establishes accessibility for systems without periodic points when the ambient 3-manifold has non-virtually solvable fundamental group. These results do not assume volume preservation, transitivity, or that every point is nonwandering. 

The classification of three-dimensional partially hyperbolic dynamics provides a route from these special classes to the general problem. Fenley and Potrie~\cite{FP26classification} proved that a chain-transitive partially hyperbolic diffeomorphism on a closed three-manifold with fundamental group of exponential growth is a collapsed Anosov flow. They have announced joint work with Hammerlindl extending the classification without chain transitivity, and hence including nonconservative systems with wandering points. This makes accessibility for general collapsed Anosov flows, without recurrence assumptions, a natural remaining question. It is the question addressed here.

\subsection{Accessibility and collapsed Anosov flows}

Removing volume preservation is useful beyond the ergodicity problem. Accessibility also allows one to compare a function along stable and unstable leaves and to determine whether these comparisons are compatible around closed $su$-paths. In Liv\v{s}ic theory, this compatibility replaces periodic-orbit obstructions that may be insufficient in partially hyperbolic dynamics. More precisely, Wilkinson~\cite{Wil13} proved that, for a $C^1$ accessible partially hyperbolic diffeomorphism and a H\"older function $a:M\to\R$, the equation
\[
a=u\circ f-u+c
\]
has a continuous solution $u$ for some constant $c\in\R$ if and only if the periodic-cycle functional vanishes on every $su$-cycle. Every continuous solution is then H\"older. Neither assertion requires volume preservation or a nonwandering hypothesis. An accessibility theorem in the nonconservative category therefore supplies a geometric hypothesis for cohomological equations and rigidity questions, as well as extending the scope of the ergodicity program.

Accessibility also plays a central role in the study of invariant measures and rigidity. For smooth cocycles with stable and unstable holonomies over $C^2$ volume-preserving center-bunched systems, the invariance principle of Avila, Santamaria and Viana~\cite{ASV13} relates vanishing fiber exponents to holonomy invariance. For cocycle-invariant measures projecting to the base volume, accessibility makes the resulting family of conditional measures continuous over the whole base. This mechanism enters the rigidity theorem of Avila, Viana and Wilkinson~\cite{AVW15} for smooth volume-preserving perturbations, sufficiently small in the $C^1$ topology, of geodesic-flow time-one maps on closed negatively curved surfaces. They prove that the disintegration of volume along center leaves is either atomic or Lebesgue, and in the latter case the perturbation is the time-one map of a smooth volume-preserving flow. Even though these results rely on volume preservation, they offer insights into exploring rigidity phenomena for more general invariant measures via accessibility.

Accessibility also helps control measures of maximal entropy. For accessible $C^{1+\alpha}$ dynamically coherent partially hyperbolic diffeomorphisms of closed three-manifolds with compact one-dimensional center leaves, Hertz, Hertz, Tahzibi and Ures~\cite{HHTU12} proved finiteness of the ergodic measures of maximal entropy, with a dichotomy between a unique measure with zero center exponent and finitely many hyperbolic ergodic maximizing measures. Crovisier and Poletti~\cite{CP} extended the invariance principle to partially hyperbolic systems with quasi-isometric center action, including systems with noncompact center leaves, and obtained applications to maximal-entropy measures and rigidity under additional hypotheses, where accessibility plays a significant role.

Collapsed Anosov flows were introduced by Barthelm\'e, Fenley and Potrie~\cite{BFP23} to describe partially hyperbolic diffeomorphisms through the orbit structure of Anosov flows. The defining data consist of a topological Anosov flow $\phi$ on $M$, a self orbit equivalence $\beta$ of $\phi$, and a continuous map $h:M\to M$ homotopic to the identity such that
\[
f\circ h=h\circ\beta.
\]
Here $\beta$ is a homeomorphism sending flow orbits to flow orbits and preserving their orientations. The map $h$ is differentiable along these orbits, with nonzero derivative tangent to $E^c$. We use the general definition in~\cite[Definition 2.7]{BFP23}. In particular, we do not require $h$ to send weak stable and weak unstable leaves to center-stable and center-unstable surfaces, as in the strong collapsed Anosov class, and we do not assume dynamical coherence.

The following theorem is the main result of this paper. Together with the announced extension of the classification by Hammerlindl, Fenley and Potrie beyond the chain-transitive setting described in~\cite{FP26classification}, it would settle Conjecture~\ref{conj:hhu-accessibility} in full generality.

\begin{theorem}\label{thm:main}
Let $f:M\to M$ be a $C^1$ collapsed Anosov flow on a closed connected three-manifold. If $\pi_1(M)$ is not virtually solvable, then $f$ is accessible.
\end{theorem}

Fenley and Potrie~\cite{FP24} proved this result under the assumption that the nonwandering set is the whole manifold. The theorem above removes that assumption for the same general class of collapsed Anosov flows, proving the nonconservative accessibility for this class. 

Without recurrence, the geometry of the nonaccessibility lamination becomes harder to control. In dimension three, nonaccessibility gives a compact lamination tangent to $E^s\oplus E^u$, but wandering regions prevent the usual recurrence arguments from controlling its complementary components. Moreover, a general collapse does not supply the weak-surface geometry available in the strong class. Our approach is to pull the lamination back by $h$ and obtain the necessary global control in the model flow. Related work of the authors~\cite{FU26a,FU26b} establishes accessibility in the absence of periodic points, including the nonconservative setting in~\cite{FU26b}. The present theorem allows periodic points and instead uses the collapsed Anosov structure.

For a three-manifold admitting a topological Anosov flow, the fundamental group is virtually solvable if and only if the flow is orbit equivalent to a suspension Anosov flow~\cite[Theorem 5.7 and Proposition 2.7]{BF13}. The hypothesis of Theorem~\ref{thm:main} can therefore be expressed by requiring that the model flow not be orbit equivalent to a suspension. Some restriction is necessary for a theorem covering the whole class. The time-one map of a constant-roof toral Anosov suspension is not accessible, since its strong stable and unstable leaves remain in the torus fibers.

This observation also gives a dynamical consequence. A \emph{discretized Anosov flow} is a partially hyperbolic diffeomorphism of the form
\[
f(x)=\psi_{\tau(x)}(x),
\]
where $\psi$ is a topological Anosov flow and $\tau:M\to\R$ is continuous~\cite[Definition 5.24]{BFP23}; see also~\cite{Martinchich} for equivalent formulations.

\begin{corollary}\label{cor:nondaf}
Every $C^1$ collapsed Anosov flow on a closed connected three-manifold for which no positive iterate is a discretized Anosov flow is accessible.
\end{corollary}

Indeed, we prove that a suspension model forces a positive iterate of $f$ to be a discretized Anosov flow. The hypothesis of Corollary~\ref{cor:nondaf} therefore excludes suspension models, and Theorem~\ref{thm:main} applies. Discretized Anosov flows with nonsuspension models are also covered by the main theorem.

\subsection{Uniform intersection with the pullback lamination}

The main geometric step is to show that a model orbit cannot avoid the pullback of a nonempty invariant $su$-lamination. Here an \emph{$su$-lamination} is a compact surface lamination with $C^1$ leaves tangent to $E^s\oplus E^u$ and saturated by both strong foliations. It may be a foliation of the whole manifold. The following statement is independent of the restriction on the fundamental group in Theorem~\ref{thm:main}.

\begin{theorem}\label{thm:intersection}
Let $f$ be a $C^1$ collapsed Anosov flow on a closed connected three-manifold, with data $(\phi,h,\beta)$. If $\Gamma$ is a nonempty invariant $su$-lamination, then there exists $T_0>0$ such that, for every $x\in M$,
\[
\phi_{(0,T_0)}(x)\cap h^{-1}(\Gamma)\ne\varnothing,
\qquad
\phi_{(-T_0,0)}(x)\cap h^{-1}(\Gamma)\ne\varnothing,
\]
where $\phi_{(0,T_0)}(x)$ denotes the flow segment between $x$ and $\phi_{T_0}(x)$.
\end{theorem}

Compact leaves of $\Gamma$ are allowed. The inverse image $h^{-1}(\Gamma)$ is a topological lamination transverse to $\phi$, and the theorem says that every open orbit segment of length $T_0$ meets this lamination. The resulting bounds on complementary flow intervals pair boundary leaves along flow segments and allow us to complete the pullback lamination to a foliation. Neither strong collapse nor recurrence is needed.

The proof also requires an obstruction to invariant circle images inside an $su$-leaf. An essential circle in the model need not have an embedded image under $h$. Proposition~\ref{prop:circle-factor} is formulated for intrinsically continuous circle factors and applies to the images arising in our construction. The intrinsic topology matters here because an $su$-leaf need not be properly embedded in $M$.

\subsection{Idea of the proof}

Suppose that $f$ is not accessible, and let $\Gamma$ be the lamination of its non-open accessibility classes. The first step is to transfer the problem to the model flow. Since $h$ has nonzero derivative in the flow direction, its inverse image $\Lambda=h^{-1}(\Gamma)$ is a transverse lamination. Each source plaque maps continuously into a specified target plaque, and $\beta$ permutes the source leaves. This local description avoids any need to classify the fibers of $h$.

To prove Theorem~\ref{thm:intersection}, assume that some model orbit misses $\Lambda$. The avoiding orbits form a compact locally maximal invariant set. Proposition~\ref{prop:freeside} produces a periodic orbit with a whole free weak side, and an iterate of $\beta$ preserves that orbit and side. Every orbit in the side meets $\Lambda$, while convergence to the periodic axis bounds the intersections at one end. The last intersection on each orbit then traces an essential circle preserved by the iterate. Its image under $h$ is intrinsically compact and lies in a single $su$-leaf. The circle-factor obstruction makes this image a point, contradicting the fact that $h$ is homotopic to the identity. This yields the uniform intersection bounds without a recurrence assumption on $f$.

The topology of $M$ excludes compact leaves of $\Gamma$. An incompressibility argument also excludes compact leaves of $\Lambda$, without requiring a source leaf to map onto an entire target leaf. The intersection bounds allow us to fill the flow gaps of $\Lambda$ and obtain a foliation without compact leaves. Its lifted leaves are proper planes, so the obstruction to accessibility can now be studied through their position relative to the weak foliations of the lifted model flow.

The next step uses a periodic lozenge in the orbit space. The opposite deck displacements of its corner orbits determine a canonical frontier plane of a union of translated slabs. The orbit projection of this plane has complete weak boundary lines, which select a deck-invariant weak intersection line. Taking the quotient gives another essential circle. When a lift of an iterate of $\beta$ preserves the lozenge, the construction is equivariant and the circle-factor obstruction again gives a contradiction.

The alternatives for self orbit equivalences in~\cite{FP24} leave two cases in which such an invariant lozenge is not supplied. In those cases, the strong-collapse classification of Fenley and Potrie~\cite{FPP24} provides suitable model data. The essential source circle then maps into one strong stable or unstable leaf, which is a line, and gives the contradiction directly. Strong collapse is used only for these alternatives, rather than asserted for an arbitrary collapsed Anosov flow.

The structure of this paper is as follows. Section~\ref{sec:prelim} fixes the definitions and geometric conventions. Section~\ref{sec:circle} proves the circle-factor obstruction and its intrinsic hyperbolic version. Sections~\ref{sec:pullback}--\ref{sec:intersection} construct the pullback lamination, establish periodic free sides, and prove uniform intersection and source properness. Section~\ref{sec:lozenge} constructs the canonical essential circles. The classification alternatives and the suspension reduction complete the proof in Section~\ref{sec:main}.

\section*{Disclosure of AI use}
The authors used AI tools to assist with editing and revising the manuscript. All ideas are attributable to the authors. The authors take full responsibility for the content of the final manuscript.

\section{Preliminaries}\label{sec:prelim}

\subsection{Partial hyperbolicity and accessibility classes}

All manifolds are closed and connected unless otherwise stated.
A $C^1$ diffeomorphism $f:M\to M$ is \emph{partially hyperbolic} if
there are a continuous $Df$-invariant splitting
\[
TM=E^s\oplus E^c\oplus E^u
\]
into three nonzero subbundles, a Riemannian metric on $M$, and an
integer $N\ge1$ such that, for every $x\in M$ and every choice of unit
vectors $v^\sigma\in E^\sigma(x)$, $\sigma\in\{s,c,u\}$,
\[
\|Df_x^N v^s\|
<
\min\{1,\|Df_x^N v^c\|\},
\qquad
\max\{1,\|Df_x^N v^c\|\}
<
\|Df_x^N v^u\|.
\]
The bundles $E^s$, $E^c$, and $E^u$ are called the stable, center,
and unstable bundles, respectively. In our three-dimensional setting,
each is one-dimensional.

The bundles $E^s$ and $E^u$ are uniquely integrable. We denote their
strong stable and strong unstable foliations by $\cW^s$ and $\cW^u$.
These foliations are $f$-invariant and have $C^1$ leaves.
An \emph{$su$-path} is a finite concatenation of arcs, each contained
in a leaf of $\cW^s$ or $\cW^u$. The \emph{accessibility class}
$AC(x)$ consists of all points that can be joined to $x$ by an
$su$-path. We say that $f$ is \emph{accessible} if $AC(x)=M$ for
every $x\in M$.

By Hertz, Hertz and Ures~\cite[Proposition A.3]{HHU08}, the set
\[
\Gamma(f)=\{x\in M:AC(x)\text{ is not open}\}
\]
is compact and $f$-invariant, and its nonopen accessibility classes
form a lamination with $C^1$ leaves tangent to $E^s\oplus E^u$.
The argument applies at $C^1$ regularity and requires no recurrence
assumption; see also~\cite[Proposition 2.4]{FU26a}.
Since $f(AC(x))=AC(f(x))$, the map $f$ permutes the leaves of this
lamination. Each leaf is saturated by both strong foliations.
We allow $\Gamma(f)=M$. If $\Gamma(f)=\varnothing$, the accessibility
classes form a partition of $M$ into open sets, so connectedness of
$M$ implies that $f$ is accessible.

We use $L,F$ and related letters for leaves of transverse laminations,
and $W$ for weak stable or weak unstable leaves of the model flow.
Curves and arcs are denoted by $c,\alpha,\ell$ or their indexed variants.
The letter $\gamma$ denotes a deck transformation, and $C$ denotes a
compact locally maximal invariant set, including the set of model orbits
that miss the pullback lamination. Stable and unstable coordinate
directions are indexed by $s$ and $u$.

A map into a leaf is \emph{intrinsically continuous} if it is continuous for
the intrinsic topology of the leaf. An intrinsically compact subset of a leaf has compact image in $M$, but the converse need not hold.

The lemma below follows from the stable manifold theorem, which produces invariant stable and unstable foliations.

\begin{lemma}\label{lem:strong-lines}
Every strong stable or strong unstable leaf is a line.
No nontrivial deck transformation preserves a lifted strong leaf.
A loop continuous in the intrinsic topology of a strong leaf is
null-homotopic in $M$.
\end{lemma}

\subsection{Topological Anosov flows}

An increasing reparametrization means
an increasing homeomorphism $\rho:\R\to\R$. We use the following
definition of a \emph{topological Anosov flow}, with the conventions
of~\cite[Definition 2.1]{BFP23}. It is a continuous flow
$\phi_t:M\to M$ satisfying the following conditions.
\begin{enumerate}
\item[(i)]
The flow is generated by a continuous nowhere-vanishing vector field
$V$. In particular, its orbits are $C^1$ curves.

\item[(ii)]
There are two invariant continuous two-dimensional foliations
$\cF^{ws}$ and $\cF^{wu}$, saturated by flow orbits and topologically
transverse. Their local product charts identify their plaques with
two transverse families of coordinate planes whose intersections
are orbit segments.

\item[(iii)]
For every $x\in M$ and $y\in\cF^{ws}(x)$, there is an increasing
reparametrization $\rho$ such that
\[
d\bigl(\phi_t(x),\phi_{\rho(t)}(y)\bigr)\longrightarrow 0
\qquad\text{as }t\longrightarrow+\infty.
\]
The corresponding statement holds for $y\in\cF^{wu}(x)$ as
$t\longrightarrow-\infty$.

\item[(iv)]
There is $\varepsilon>0$ such that, whenever
$y\in\cF^{ws}_{\varepsilon}(x)$ lies outside the orbit of $x$,
every increasing reparametrization $\rho$ satisfies
\[
d\bigl(\phi_t(x),\phi_{\rho(t)}(y)\bigr)>\varepsilon
\]
for some $t\leq 0$. Here $\cF^{ws}_{\varepsilon}(x)$ denotes a local
weak stable plaque. For
$y\in\cF^{wu}_{\varepsilon}(x)$ outside the orbit of $x$, the
corresponding separation occurs for some $t\geq 0$.
\end{enumerate}
The local product structure and the separation condition imply
expansiveness. The flow is nonsingular, and its weak foliations
have no singular prongs. 

Let us mention that, if a topological Anosov flow in dimension three is transitive, then it is orbit equivalent to a smooth Anosov flow \cite{Shannon}. Throughout this paper, the topological Anosov flows we are considering are, a priori, not transitive or smooth.

The orbit-space results of Fenley~\cite{Fen94} and
Barbot~\cite{Bar95}, together with their pseudo-Anosov extension
in Fenley--Mosher~\cite{FM01}, give the topology used below.
The universal cover $\wt M$ is homeomorphic to $\R^3$, and the
lifted flow is topologically a product. Its orbit space
\[
\cO=\wt M/\wt\phi
\]
is homeomorphic to $\R^2$. We write $q:\wt M\to\cO$ for the quotient
map and $\mathcal G^s,\mathcal G^u$ for the induced one-dimensional
foliations. Lifted weak leaves are properly embedded planes, and
their images under $q$ are properly embedded lines. Thus the orbit
quotient of a fixed weak plane is a line, equipped with either of
its two orders. A lifted weak stable plane and a lifted weak
unstable plane intersect in at most one orbit.

Weak leaves downstairs are incompressible~\cite{FM01}. After
orienting the weak directions, they are planes or cylinders, and
each cylindrical weak leaf contains a unique periodic orbit.
For smooth Anosov flows, this is the weak-leaf classification
of~\cite[Proposition 1.2]{Ver74}. The same argument applies here
using contracting weak holonomy. Uniqueness of the periodic orbit
also follows directly from the asymptotic property within a weak
leaf.

A \emph{positive self orbit equivalence} of $\phi$ is a
homeomorphism $\beta:M\to M$ carrying each flow orbit onto a flow
orbit and preserving its orientation. Equivalently, there is a
continuous function $\tau:M\times\R\to\R$ such that
\[
\beta\bigl(\phi_t(x)\bigr)
 =\phi_{\tau(x,t)}\bigl(\beta(x)\bigr),
\qquad
\tau(x,t+s)=\tau(x,t)+\tau(\phi_t(x),s),
\]
where $\tau(x,0)=0$ and $\tau(x,\cdot)$ is an increasing
homeomorphism of $\R$ for every $x\in M$. Throughout the paper,
self orbit equivalences are understood to be positive.

Such a map preserves both weak foliations. Indeed, uniform
continuity carries asymptotic time-matched semiorbits to
asymptotic time-matched semiorbits, and the increasing time changes
preserve the relevant time direction. The characterization of
weak leaves by these asymptotic relations therefore gives
\[
\beta\bigl(\cF^{ws}(x)\bigr)
 \subset \cF^{ws}(\beta(x)),
\qquad
\beta\bigl(\cF^{wu}(x)\bigr)
 \subset \cF^{wu}(\beta(x)).
\]
Applying the same argument to $\beta^{-1}$ gives equality.

\subsection{Collapse data and finite reductions}\label{sec:CAF-data}

Following~\cite[Definition 2.7]{BFP23}, we call $f$ a
\emph{collapsed Anosov flow} if there are a topological Anosov
flow $\phi$, a continuous map $h:M\to M$ homotopic to the
identity, and a positive self orbit equivalence $\beta$ of
$\phi$ such that
\[
f\circ h=h\circ\beta
\]
and, for every $x\in M$, the map $t\mapsto h(\phi_t(x))$ is
differentiable with
\[
\left.\frac{d}{dt}\right|_{t=0}h(\phi_t(x))
 \in E^c(h(x))\setminus\{0\}.
\]
We call $(\phi,h,\beta)$ the collapse data. Continuity of the
orbitwise derivative is not required.

Since $h\simeq\id$, the map $h$ is onto. Its lift obtained from
an identity homotopy,
\[
\wt h:\wt M\longrightarrow \wt M,
\]
commutes with every deck transformation and stays a bounded
distance from the identity. Consequently, $\wt h$ is proper.
No injectivity or connectedness of its fibers will be used.

Orbitwise differentiability gives the geometric regularity needed
for pullback charts. Locally choose a continuous nonzero oriented
center field $Y$ and a coordinate with $Y_j\ne0$. An orbit image
$c$ satisfies
\[
c'(t)=b(t)Y(c(t)),
\qquad b(t)\ne0.
\]
The derivative $c_j'(t)$ has a constant sign by the
intermediate-value property of derivatives. Thus $c_j$ is
strictly monotone, and reparametrizing by $\zeta=c_j(t)$ gives
\[
\frac{dc_i}{d\zeta}=\frac{Y_i(c)}{Y_j(c)}.
\]
The geometric center arc is therefore $C^1$, even when the
original orbit speed is not assumed continuous. Its transversality
to the $su$-plaques is uniform on compact coordinate charts.

When required, we pass to a common finite cover orienting the
manifold, the strong bundles, and the model weak transverse
directions. This cover may be chosen characteristic. Indeed, a
finitely generated group has only finitely many subgroups of a
given finite index, so the intersection of all automorphic images
of a finite-index subgroup is a characteristic finite-index
subgroup contained in it.

The maps $f$ and $\beta$ have the same outer action on
$\pi_1(M)$ by $f\circ h=h\circ\beta$ and $h\simeq\id$.
Lift $h$ along its identity homotopy and choose compatible lifts
of $f$ and $\beta$ so that the collapse identity remains exact.
A positive iterate preserves the relevant orientations. The
same construction works in any prescribed finite cover after
taking a suitable iterate, since the induced outer automorphism
permutes finitely many conjugacy classes of subgroups of that
index.

Positive iterates have the same strong foliations. Accessibility
of a lifted iterate descends by projecting strong paths, and
non-virtual solvability is preserved in finite-index subgroups.
These reductions consequently preserve the contradiction argument
for Theorem~\ref{thm:main}. They also preserve the intrinsic leaf
maps and semiconjugacies constructed below.

We will use strong collapse only at the end of the proof.
Following~\cite[Definition 2.10]{BFP23}, a collapsed Anosov flow
is a \emph{strong collapsed Anosov flow} if its collapse data
can be chosen to satisfy two additional requirements.
\begin{enumerate}
\item[(i)]
For every lifted weak stable leaf $W$, the image $\wt h(W)$ is a
properly embedded $C^1$ plane tangent to $E^{cs}$. For a lifted
weak unstable leaf, the corresponding image is tangent to
$E^{cu}$.

\item[(ii)]
For every lifted flow orbit $o$ contained in such a weak plane $W$, the curve
$c=\wt h(o)$ separates $\wt h(W)$ into exactly two components. Each
component of $W\setminus o$ is mapped into the closure, taken
inside $\wt h(W)$, of exactly one component of $\wt h(W)\setminus c$.
\end{enumerate}
The second condition is the transverse collapsing condition.
It concerns the two sides of an orbit within its weak leaf.

In particular, a source curve contained in a weak plane and in
one pullback $su$-leaf has its image locally in the intersection
of two specified transverse $C^1$ surfaces. This intersection
is tangent to $E^s$ in the weak stable case and to $E^u$ in the
weak unstable case. This is the additional strong-collapse
property needed below.

\begin{lemma}\label{lem:no-compact}
Let $f$ be a collapsed Anosov flow on a closed three-manifold
with non-virtually-solvable fundamental group. No
$su$-lamination for $f$ has a compact leaf.
\end{lemma}

\begin{proof}
Pass to a finite cover orienting the manifold and the strong
bundles, and take a compatible lifted iterate. A compact
$su$-leaf is then an orientable closed surface carrying a
nowhere-vanishing stable vector field, and hence is a torus.

By~\cite[Proposition 2.4]{HHU16}, the existence of an $su$-torus
implies that the covering manifold admits an Anosov torus.
The covering manifold is irreducible because its universal
cover is $\R^3$. The classification in~\cite[Theorem 1.1]{HHU11} therefore gives the three-torus, the mapping torus of $-\id$ on the two-torus, or the mapping torus of a hyperbolic toral automorphism. Each has virtually solvable fundamental group, contradicting the preservation of non-virtual solvability under passage to finite covers.
\end{proof}

\section{Expansive factors of circle homeomorphisms}\label{sec:circle}

A homeomorphism $G$ of a compact metric space $(K,d)$ is \emph{expansive} if there is $\delta>0$ such that if $d(G^nx,G^ny)<\delta$ holds for every $n\in\Z$, then we have $x=y$. Let us first introduce the following result about circle homeomorphisms with expansive factors. The factor map need not have connected fibers. The argument first collapses intervals with shrinking two-sided iterates, and then treats the remaining rotation factor.

\begin{proposition}\label{prop:circle-factor}
Let $g:S^1\to S^1$ be a homeomorphism and let $G:K\to K$ be an expansive homeomorphism of a compact metric space. If a continuous surjection $\eta:S^1\to K$ satisfies $G\circ \eta=\eta\circ g$, then $K$ consists of one point.
\end{proposition}

\begin{proof}
Fix an expansivity constant $\delta>0$ for $G$. Every positive power $G^r$ is expansive, since uniform continuity of $G^j$, $0\le j<r$, converts a sufficiently small expansivity scale for the times $rn$ into the scale $\delta$ for all integer times. We may therefore replace $g,G$ by common positive powers, and first take a square so that $g$ preserves orientation.

Suppose that a compact interval $J_0\subset S^1$ satisfies
\[
\operatorname{diam}(g^nJ_0)\longrightarrow0
\qquad\text{as }n\longrightarrow\pm\infty.
\]
Then $\eta$ is locally constant on $\Int J_0$. To see this, choose $\delta_0>0$ such that circle points at distance less than $\delta_0$ have $\eta$-images at distance less than $\delta$. All sufficiently large positive and negative iterates of $J_0$ have diameter less than $\delta_0$. Around any interior point, choose an interval $J\subset J_0$ whose remaining finitely many iterates also have diameter less than $\delta_0$. For $s,t\in J$ we have
\[
d(G^n\eta(s),G^n\eta(t))
=d(\eta(g^ns),\eta(g^nt))<\delta
\qquad(n\in\Z).
\]
Expansivity gives $\eta(s)=\eta(t)$.

Assume first that $g$ has a periodic point. Taking a positive power, let $P=\operatorname{Fix}(g)\ne\varnothing$. Every component $I$ of $S^1\setminus P$ is an invariant interval on which $g$ is increasing and fixed-point free. The iterates of a compact subinterval $I_0\subset I$ approach one endpoint in positive time and the other in negative time. They satisfy the property $\operatorname{diam}(g^nI_0)\longrightarrow0$ as $n\to \pm\infty$, so $\eta$ is locally constant, hence constant, on $I$. Continuity gives the same value at its endpoints.

For $x\in P$, semiconjugacy gives $\eta(x)\in\operatorname{Fix}(G)$. The latter set is finite, since distinct fixed points are $\delta$-separated. Therefore
\[
\eta(S^1)=\eta(P)\subset\operatorname{Fix}(G).
\]
The image of the circle is connected, so it is a point.

Assume now that $g$ has no periodic point. By Poincar\'e's classical rotation theory, $g$ has irrational rotation number and is semi-conjugate to an irrational rotation $R_\theta: S^1\to S^1$, with $\theta\in\R\setminus\mathbb Q$, through a monotone surjective map $\varpi:S^1\to S^1$.

For any point $z\in S^1$, it is clear that $\eta(\varpi^{-1}(z))\in K$ is a single point if $\varpi^{-1}(z)$ is a single point. If $I=\varpi^{-1}(z)$ is a nontrivial fiber, its iterates $g^nI$ are distinct compact arcs with pairwise disjoint interiors. Their lengths tend to zero as $|n|\to\infty$. The preceding interval argument and continuity show that $\eta$ is constant on $I$. Thus, the composition $\psi=\eta\circ\varpi^{-1}:S^1\to K$ gives a well-defined continuous surjection such that
\[
\eta=\psi\varpi,
\qquad G\psi=\psi R_\theta.
\]
Uniform continuity supplies $\vartheta_0>0$ such that
\[
d(\psi(z+\vartheta),\psi(z))<\delta
\qquad(z\in S^1,\ |\vartheta|<\vartheta_0).
\]
Consequently, for every $n\in\Z$,
\[
d(G^n\psi(z+\vartheta),G^n\psi(z))
=d(\psi(z+n\theta+\vartheta),\psi(z+n\theta))<\delta.
\]
Expansivity makes $\psi$ locally constant. It is constant by connectedness, and its surjectivity proves that $K$ is a point.
\end{proof}

The application concerns hyperbolicity in the intrinsic topology of an $su$-leaf. Ambient properness of the leaf is unnecessary.

The following lemma is a well-known fact using local product structure for stable and unstable manifolds. We do not present its proof here.

\begin{lemma}\label{lem:intrinsic-expansivity}
Let $f$ be a $C^1$ partially hyperbolic diffeomorphism. Let $L$ be a $C^1$ surface tangent to $E^s\oplus E^u$, saturated by both strong foliations, and suppose $f^m(L)=L$ for some $m>0$. If $K\subset L$ is intrinsically compact and $f^m(K)=K$, then $f^m|_K$ is expansive for the intrinsic metric of $L$.
\end{lemma}

\begin{corollary}\label{cor:su-circle}
Under the hypotheses on $f,L,m$ in Lemma~\ref{lem:intrinsic-expansivity}, let $g:S^1\to S^1$ be a homeomorphism. Every continuous map $\eta:S^1\to L$ satisfying $f^m\circ\eta=\eta\circ g$ is constant.
\end{corollary}
\begin{proof}
The set $K=\eta(S^1)$ is intrinsically compact, and $f^m(K)=K$ because $g$ is onto. Apply Lemma~\ref{lem:intrinsic-expansivity} and Proposition~\ref{prop:circle-factor} to $f^m|_K$.
\end{proof}

\section{Pulling back the accessibility lamination}\label{sec:pullback}

Fix collapsed Anosov flow data $(\phi,h,\beta)$ for $f$. Let $\Gamma$
be a nonempty compact invariant $su$-lamination, with its leaves
permuted by $f$, and put
\[
 \Lambda=h^{-1}(\Gamma).
\]
We allow $\Gamma=M$. Intrinsic continuity always implies ambient continuity.
The converse holds when the domain is locally connected, but may
fail for a general domain when a leaf accumulates on itself.
In particular, the two notions are equivalent for maps
between leaves when the source leaf carries its intrinsic
manifold topology.

\begin{proposition}\label{prop:pullback}
The set $\Lambda$ is a nonempty compact topological lamination transverse
to $\phi$. The restriction of $h$ to each source leaf is intrinsically
continuous into one leaf of $\Gamma$. Moreover, $\beta$ permutes the
leaves of $\Lambda$ by intrinsic homeomorphisms.
\end{proposition}

\begin{proof}
Since $h\simeq\id$, the map $h$ is onto. Thus $\Lambda$ is nonempty
and compact. The collapse identity and invariance of $\Gamma$ give
\[
 \beta(\Lambda)=\Lambda.
\]
We verify the local lamination structure and then its leafwise
invariance.

Choose a target lamination box whose plaques are ordered $C^1$ graphs
over a disk, transverse to the center direction. By the discussion in Section~\ref{sec:CAF-data}, a sufficiently short image of a model
orbit is a regular center arc geometrically. In a smaller box its
tangent stays in a fixed cone transverse to the plaque graphs. It
therefore crosses each plaque at most once, in the prescribed order.
This uses only the geometric center direction; no continuity of the
original orbit speed is needed.

Take a source flow box around $x\in\Lambda$ small enough that its
image is in this target box. Two short flow times on either side of
$x$ bracket the crossing of the target plaque through $h(x)$.
Continuity of $h$ preserves the bracketing for nearby source fibers
and nearby target plaque parameters. The crossing time is unique and
continuous in both parameters. Indeed, for convergent parameters the
times stay in the fixed bracketing interval; every convergent
subsequence has, by continuity of $h$, the prescribed crossing as its
limit, and uniqueness identifies that limit. The ordered crossing
graphs are exactly the local plaques of $\Lambda$.

On one such graph the target plaque parameter is fixed. The target
disk coordinate of $h$ is continuous, so its restriction is continuous
into that specified plaque in the intrinsic topology. A path in a
source leaf is covered by finitely many pullback boxes. On overlaps
the specified target plaques belong to the same leaf. Continuing along
the path proves that $h$ maps the whole source leaf intrinsically
continuously into one target leaf. The same plaque argument holds for
the lifted maps and leaves in the universal cover.

It remains to establish leafwise invariance. Let $P$ be a small source
plaque whose image lies in a specified target plaque $P_{\Gamma}$. The
restriction of $f$ between target leaves is intrinsically continuous.
Shrink $P$ so that $f(h(P))$ lies in one target plaque $P_{\Gamma}'$, and then
so that $\beta(P)$ is contained in one source pullback box. Since
\[
 h\beta(P)=fh(P)\subset P_{\Gamma}',
\]
the set $\beta(P)$ lies in the unique crossing graph belonging to
$P_{\Gamma}'$. The map into this source plaque is intrinsically continuous.
The same argument applied to $\beta^{-1}$, using
$h\beta^{-1}=f^{-1}h$, gives the inverse local assertion. Continuation
along source leaves proves that $\beta$ permutes them by intrinsic
homeomorphisms. 
\end{proof}

We will not require a source leaf to map onto its target leaf. Moreover, we do not assume connectedness or monotonicity of the fibers of $h$. Intrinsic continuity is enough to exclude compact source leaves when the target leaves are noncompact.

\begin{lemma}\label{lem:no-compact-source}
If every leaf of $\Gamma$ is noncompact, then every leaf of $\Lambda$
is noncompact.
\end{lemma}

\begin{proof}
Pass to a common finite cover orienting the ambient manifold and the model weak transverse directions, and retain the notation. The lifted target leaves remain noncompact. If the original pullback had a compact leaf, a component of its lift would be compact, so it suffices to prove the assertion in this cover.
Suppose that $T$ is a compact leaf of $\Lambda$. The model weak foliations induce nonsingular one-dimensional foliations on $T$, so $\chi(T)=0$. Since $T$ is orientable, it is a torus.

We first show that $T$ is incompressible. The classical
incompressibility result for transverse tori appears in
Mosher~\cite{Mos92} and Brunella~\cite{Bru93}; see also~\cite{Pot25}.
We give an orbit-space argument that applies directly to the
topologically transverse torus considered here.

Let $\pi:\wt M\to M$ be the universal covering map, and let $S$ be a
component of $\pi^{-1}(T)$. Since $T$ is locally a transverse graph,
$S$ is a properly embedded, locally flat surface in $\wt M\cong\R^3$.
It separates $\wt M$, and the crossing direction of $\wt\phi$ is constant
along $S$. Thus every lifted orbit meets $S$ at most once.
Consequently, the restriction of the orbit projection
\[
q|_S:S\longrightarrow\cO
\]
is an injective local homeomorphism and identifies $S$ with an open
subset $U\subset\cO\cong\R^2$.

Write $i:T\hookrightarrow M$ for the inclusion and $i_*:\pi_1(T)\to \pi_1(M)$ for the induced action. The covering $S\to T$ corresponds to $\ker i_*$. If $i_*$ were not injective, then $S$ would be a torus or an open annulus. The torus case is impossible because $q(S)$ is open in the plane.

Suppose that $S$ is an annulus. Its deck group contains an infinite
cyclic subgroup of finite index. Choose a generator $\gamma$ of this
subgroup, and use the same notation for its action on $\cO$.
Let $c\subset U$ be an essential simple closed curve. Properness
of the action on $S$ implies that $\gamma^n(c)\cap c=\varnothing$ for
all sufficiently large $n$. These two curves are essential in the
same planar annulus, so their bounded Jordan disks are nested.
After replacing $\gamma$ by a power or its inverse, the disk $D$ bounded
by $c$ therefore satisfies
\[
\gamma(D)\subset\operatorname{int}D.
\]

By the Brouwer fixed-point theorem, $\gamma$ has a fixed point
$p\in\operatorname{int}D$. The action at $p$ is a topological saddle.
Let $\ell$ be its repelling invariant leaf of $\mathcal G^s$ or
$\mathcal G^u$. The weak-leaf properties established above imply
that $p$ is the only fixed point of $\gamma$ on $\ell$. Indeed, the
corresponding weak leaf downstairs is an incompressible cylinder
with a unique periodic orbit. This orbit represents the core of
the cylinder and has a single lift in the corresponding lifted
weak plane.

Since $\ell$ is a properly embedded line, the forward iterates
of any $x\ne p$ sufficiently close to $p$ on $\ell$ escape every
compact subset of $\cO$. Choosing $x\in\operatorname{int}D$
contradicts $\gamma^k(x)\in D$ for all $k\ge0$. This excludes the annulus
case and proves that $i_*$ is injective.

By Proposition~\ref{prop:pullback}, the restriction $h_T:T\to L$
is intrinsically continuous for some leaf $L$ of $\Gamma$.
Its image is intrinsically compact and hence lies in a compact
connected subsurface $K\subset L$ with nonempty boundary.
The composition
\[
\Z^2=\pi_1(T)
\xrightarrow{(h_T)_*}\pi_1(K)
\longrightarrow\pi_1(M)
\]
agrees with $i_*$ up to basepoint conjugacy because $h\simeq\id$.
It is therefore injective, so $(h_T)_*$ is injective.
However, since $K$ is a subsurface of a non-compact leaf $L\in\Gamma$, $\pi_1(K)$ is free and contains no subgroup isomorphic
to $\Z^2$. This contradiction proves the lemma.
\end{proof}

\section{Periodic free sides}\label{sec:free}

Throughout this section, $\phi=(\phi_t)_{t\in\mathbb R}$ is a
topological Anosov flow on a closed connected three-manifold $M$.
Its weak stable and weak unstable foliations are denoted by
$\mathcal F^{ws}$ and $\mathcal F^{wu}$.

A compact invariant set $C$ is \emph{locally maximal} if
\[
C=\Inv_\phi(N):=\bigcap_{t\in\R}\phi_t(N),
\qquad C\subset\Int N,
\]
for some compact neighborhood $N$. Let $p\colon\widehat M\to M$ be a connected finite covering in which the directions transverse to the flow inside both weak
foliations are orientable. Denote by $\widehat\phi$ the lifted flow
and set
\[
  \widehat C=p^{-1}(C).
\] The weak stable leaf of any $\widehat \phi$-periodic orbit is a cylinder. An \emph{entire stable side} is a component of this cylinder with the periodic orbit removed; unstable sides are defined analogously.

The orientation condition means that the one-dimensional stable
and unstable foliations induced on transverse sections admit
orientations compatible under changes of flow box. Such a finite
covering always exists; if these orientations already exist on $M$,
one may take $p=\operatorname{id}_M$.
Here $p^{-1}(C)$ is the full set-theoretic preimage, including all
its components. 

The following proposition is our main result of this section, whose proof is deferred to the end of the section.

\begin{proposition}\label{prop:freeside}
Let $\phi$ be a topological Anosov flow on a closed connected three-manifold, and let $C$ be a nonempty proper compact locally maximal invariant set. In a finite cover $\widehat M$ orienting the weak directions, the set of periodic orbits of $\widehat \phi$ contained in the full inverse image $\widehat C=p^{-1}(C)$ that have at least one entire stable or unstable side disjoint from $\widehat C$ is nonempty and finite.
\end{proposition}

The set $C$ need not be transitive or have empty interior. In fact, for each choice of weak direction and side, the points with an entire free side form a finite union of periodic orbits. This finiteness makes the construction compatible with a self orbit equivalence preserving $C$.

For the proofs, pass to the finite covering above and replace
$C\subset M$ by its full preimage
\[
  \widehat C=p^{-1}(C)
  =\{\widehat x\in\widehat M\mid p(\widehat x)\in C\}.
\]
For the remainder of the proof we relabel
$(\widehat M,\widehat\phi,\widehat C)$ as $(M,\phi,C)$.
We write $\pi\colon\wt M\to M$ for the universal covering,
$\wt\phi$ for the lifted flow, and
$\wt C=\pi^{-1}(C)$.
Thus preimages always include all components. The notation $p^{-1}(C)$
or $\pi^{-1}(C)$ denotes a set-theoretic preimage, not an inverse map.

Each lifted weak leaf $W$ is a plane, and the quotient of $W$ by its
flow orbits is a line. The chosen transverse orientation orders this
line. If $q_W\colon W\to W/\wt\phi$ is the quotient map and
$W=\wt{\mathcal F}^{ws}(\wt x)$, set
\[
  P_s^+(\wt x)
  =q_W^{-1}\bigl(\{\xi\in W/\wt\phi
                     \mid \xi>q_W(\wt x)\}\bigr).
\]
This is the positive component of
$W\setminus\wt\phi_{\mathbb R}(\wt x)$.
Define $P_s^-(\wt x)$ by the opposite inequality, and define
$P_u^\pm(\wt x)$ in the lifted weak unstable leaf in the same way.
These components are open half-planes in the intrinsic topology of
the leaf. Deck transformations preserve the chosen orientations.

We use product coordinates $(s,u)\in J^s\times J^u$ on a small
transverse rectangle $R$. The intervals $J^s\times\{u\}$ and
$\{s\}\times J^u$ are respectively its stable and unstable plaques,
denoted by $R^s(z)$ and $R^u(z)$ when they contain $z$.
In particular, for $z=(s_z,u_z)$ and $w=(s_w,u_w)$,
\[
  [z,w]_R=R^s(z)\cap R^u(w)=(s_w,u_z).
\]
The $s$-coordinate increases in the chosen positive stable direction.
For a flow box based on $R$, the projection along short flow segments
is the map
\[
  \operatorname{pr}_R(\phi_t(z))=z,
  \qquad z\in R,\quad |t|<\delta_R,
\]
where $\delta_R>0$ is chosen so that the flow-box coordinates are unique.

We also specify the meaning of transport between transverse sections.
Given an orbit segment from $z\in R$ to $\phi_T(z)\in R'$, its
Poincar\'e transport $H\colon R\to R'$ is the continuation of local flow-box projections
along that segment. On an interval $J^{\sigma}\subset R$ where this continuation
is defined and $\sigma=s,u$, it has the form
\[
  H(w)=\phi_{\tau(w)}(w),
  \qquad \tau\colon J^{\sigma}\to\mathbb R\ \text{continuous},
  \qquad \tau(z)=T,
  \qquad H(J^{\sigma})\subset R'.
\]
The specified orbit segment determines the branch of the continuation.
The map is not generally the time-$T$ map, since $\tau(w)$ may depend
on $w$. 

Choose a lift $\wt z\in \wt J^{\sigma}$ so that $\pi|_{\wt J^{\sigma}}:\wt J^{\sigma}\to J^{\sigma}$ is a homeomorphism. We use the lifted flow segment from $\wt z$ to $\wt\phi_T(\wt z)$ and define the corresponding lifted transport by
\[
  \wt H(\wt w)
  =\wt\phi_{\tau(\pi(\wt w))}(\wt w).
\]
If this transport extends along a larger lifted arc, we cover the arc by overlapping intervals $\wt J^{\sigma}_j$ on which $\pi$ is injective and continue the same transport branch from the interval containing $\wt z$. The resulting local maps $\wt H_j$ agree on overlaps, and hence define a single transport map on the arc. Throughout this continuation, the terminal section lift is the one determined by $\wt\phi_T(\wt z)$. Even if the projected base orbit returns near an earlier point of $M$, its lifted endpoint may lie in a different lift of the same transverse section.

Fix a compatible metric $d$ on $M$.
By an increasing time change we mean an increasing homeomorphism
between the time intervals under consideration.

\begin{lemma}\label{geo:contraction}
The weak foliations have the following properties.
\begin{enumerate}
\item
For every $\eps>0$, sufficiently close points $z,w$ in a common
local weak stable plaque admit an increasing homeomorphism
$\rho\colon[0,\infty)\to[0,\infty)$ with $\rho(0)=0$ such that
\[
  d(\phi_t(z),\phi_{\rho(t)}(w))<\eps
  \qquad\text{for all }t\geq0.
\]
The required closeness can be chosen uniformly in a finite foliation
atlas. The analogous statement holds on $(-\infty,0]$ for local weak
unstable plaques.

\item
Fix a sufficiently small finite system of transverse rectangles
$R_i\subset\operatorname{int}S_i$, whose inner rectangles meet every
orbit segment of some fixed length. For $z\in R_i$ and
$z_T=\phi_T(z)\in R_j$, we denote by $R^s_i(z)$ and $S^s_j(z_T)$ the corresponding stable plaques given by the product coordinates as above. Then, for all sufficiently large positive
visit times $T$, the transport along this orbit segment is defined
on the whole interval $R_i^s(z)$ and satisfies
\[
  H_{z,T}(R_i^s(z))\subset S_j^s(z_T),
\]
and
\[
  \sup_{w\in R_i^s(z)}
     |s(H_{z,T}(w))-s(z_T)|\longrightarrow0
  \quad\text{as }T\longrightarrow\infty,
\]
where $s(\cdot)$ is the stable coordinate. The convergence is uniform in $i,j,z$ and the chosen visit time.
The corresponding assertion holds for unstable intervals as
$T\to-\infty$.

\item
Let $\alpha$ be a fixed compact arc transverse to the flow in a lifted
weak stable plane, and fix $\wt z\in \alpha$. At sufficiently late
visits of the orbit of $\wt z$ to lifts of the inner rectangles,
transport along that lifted orbit extends over $\alpha$. Its image lies
in a single stable plaque of the corresponding lifted outer rectangle,
and its stable-coordinate width tends to zero.
The threshold may depend on $\alpha$.
\end{enumerate}
The conclusions are unchanged by composing with fixed local flow-box
projections at the initial or terminal section, after restricting
their domains as necessary.

Moreover, if $C$ is compact and locally maximal, then in every
sufficiently small transverse product rectangle,
\[
  C\cap R=K^s\times K^u,
\]
where $K^s\subset J^s$ and $K^u\subset J^u$ are relatively closed.
\end{lemma}

\begin{proof}

The local stability and uniform contraction statements follow from
\cite[Lemmas 5--7]{Iak25}. Since the time gaps between successive
visits to the inner rectangles are uniformly bounded, the number
of visits tends uniformly to infinity with the elapsed transport
time. Thus the contraction estimates stated in terms of the number
of visits also apply for all sufficiently large transport times.
The outer rectangles provide the neighborhoods needed to continue
the chosen Poincar\'e branches near the boundaries of the inner
rectangles. The transported interval contracts inside the terminal
local stable plaque; no return to the initial rectangle is required.
Moreover, the product-coordinate maps are uniformly continuous
on the finitely many compact rectangles under consideration.
Consequently, the stable-coordinate diameter of the transported
interval also tends uniformly to zero.

To pass from transverse sections to local weak plaques, we project
along short flow segments inside flow boxes. These projections
introduce small time offsets, and any resulting short initial orbit
segments are controlled by uniform continuity of the flow.
For the local stability assertion, a matching time change
$\rho_\delta:[0,\infty)\to[\delta,\infty)$ with $\rho_\delta(0)=\delta$ can be
normalized by setting
\[
  \rho(t)=\rho_\delta(t)-\delta.
\]
The map $\rho$ is an increasing homeomorphism of $[0,\infty)$ fixing zero, and
\[
  d\bigl(\phi_{\rho_\delta(t)}(y),\phi_{\rho(t)}(y)\bigr)
  \le \sup_{z\in M}d\bigl(\phi_\delta(z),z\bigr)
  \longrightarrow 0
  \qquad\text{as }\delta\longrightarrow0,
\]
uniformly for $t\ge0$, by compactness of $M$ and continuity
of the flow. Choosing the initial matching error and the flow
offsets sufficiently small therefore gives the required stability
estimate with a time change fixing zero. Finally, the additional
local flow-box projections preserve convergence by uniform
continuity on the compact flow boxes used in the construction.

For the third assertion, cover $\alpha$ by a finite chain of overlapping
transverse stable intervals small enough for local stability to apply.
Start with the interval containing $\wt z$. Once it has been
transported, an overlap selects the continuation on the next interval.
The elapsed time at that overlap tends to infinity with the base
time, because the stable matching time change is increasing and onto.
Local contraction therefore applies on the next interval as well.
Induction over the finite chain, using terminal widths small enough
for their union to remain inside one outer rectangle, gives the
assertion. The maps agree on overlaps because they continue the same
local flow holonomy. This is the finite-chain argument used in
\cite[Lemma 7]{Iak25}, now applied to a fixed arc rather than to a
uniformly bounded family of local intervals.

For the final assertion, take a compact neighborhood $N$ such that
\[
  C=\operatorname{Inv}_{\phi}(N)
   :=\{z\in N\mid \phi_t(z)\in N\text{ for every }t\in\mathbb R\}
   \subset\operatorname{int}N.
\]
Choose $\eps>0$ so that the $\eps$-neighborhood of $C$ lies
in $N$. If $z,w\in C$ lie in a sufficiently small product rectangle,
the positive semiorbit of $[z,w]_R$ stays within $\eps$ of the
orbit of $z$, and its negative semiorbit stays within $\eps$ of
the orbit of $w$, with the time changes in the first assertion.
Thus its entire orbit lies in $N$, and $[z,w]_R\in C$.
Closure under both coordinate combinations gives the product identity.
The factors are relatively closed because $C$ is closed.
\end{proof}

The uniformity in the second assertion is needed for local intervals
that vary along a sequence. The third assertion applies to a global
arc only after that arc has been fixed.

For a compact locally maximal invariant set $C$, define
\[
  \mathscr B_s^+(C)=\{x\in C\mid
       P_s^+(\wt x)\cap\wt C=\varnothing\},
\]
where $\wt x$ is any lift of $x$.
This definition is independent of the lift, since deck transformations
preserve both $\wt C$ and the positive stable direction.
It means that the entire positive side in the lifted weak stable leaf,
rather than merely a local interval, is disjoint from $\wt C$.

\begin{lemma}\label{geo:free-points}
The set $\mathscr B_s^+(C)$ is compact and invariant. If it is nonempty, it is
a finite union of periodic orbits. The same conclusion holds for
$\mathscr B_s^-(C)$ and $\mathscr B_u^\pm(C)$, defined using the other choices of side
and weak foliation.
\end{lemma}

\begin{proof}
Invariance follows because the defining condition depends only on
the lifted orbit. To prove closedness, let $x_n\in \mathscr B_s^+(C)$ converge
to $x\in C$. If $x\notin \mathscr B_s^+(C)$, choose
$\wt w\in P_s^+(\wt x)\cap\wt C$ and a compact
transverse arc starting at $\wt x$ and ending on the orbit of
$\wt w$ in their common weak stable plane.
By Lemma~\ref{geo:contraction}(3), a fixed sufficiently late
transport puts the endpoints in the interior of one lifted product
rectangle. Their projections have coordinates $x'=(s_0,u_0)$ and $w'=(s_1,u_0)$ with $s_0<s_1$.
Choose the rectangle small enough for the product identity for $C$.
Project $x_n$ along short flow segments to the initial section through $x$ and
then apply the same fixed local transport branch near $x$.
This gives $x_n'=(s_n,u_n)\in C$ and $x_n'\longrightarrow x'$. For large $n$, the product identity gives
\[
  z_n=(s_1,u_n)\in C,\qquad s_n<s_1.
\]
Lifting the local stable interval from $x_n'$ to $z_n$ puts its
second endpoint in $P_s^+(\wt x_n')\cap\wt C$.
This contradicts $x_n'\in \mathscr B_s^+(C)$, which follows from invariance.
Thus $\mathscr B_s^+(C)$ is closed, and hence compact.

In a sufficiently small product rectangle, all points of
$\mathscr B_s^+(C)$ have the same stable coordinate. Indeed, if two such
coordinates were $s_-<s_+$, combining $s_+$ with the unstable coordinate
of the first point would give a point of $C$ on its positive stable
side. Thus $\mathscr B_s^+(C)\cap R$ lies in one unstable plaque.
By invariance, its trace in the associated flow box lies in one
local weak unstable plaque.

We next verify positive expansivity of the restricted flow.
Fix a time tolerance $\varepsilon_0>0$. Expansivity of $\phi$ gives $e>0$ such that, for any
increasing homeomorphism $\rho\colon\mathbb R\to\mathbb R$ fixing zero, if $d(\phi_t(v),\phi_{\rho(t)}(w))<e$ holds for all $t\in\mathbb R$, then $w=\phi_\tau(v)$ for some $|\tau|<\varepsilon_0$. Compactness, the preceding one-plaque property, and
Lemma~\ref{geo:contraction}(1) give $\delta\in(0,e)$ such that
any $v,w\in \mathscr B_s^+(C)$ with $d(v,w)<\delta$ admit an increasing
homeomorphism $\rho_-\colon(-\infty,0]\to(-\infty,0]$ fixing zero with
\[
  d(\phi_t(v),\phi_{\rho_-(t)}(w))<e,\qquad\forall t\leq0.
\]
If an increasing homeomorphism
$\rho_+\colon[0,\infty)\to[0,\infty)$ fixing zero also satisfies $d(\phi_t(v),\phi_{\rho_+(t)}(w))<\delta$ for all $t\geq0$, then the two time changes paste at zero to give a homeomorphism of
$\mathbb R$. Expansivity yields $w=\phi_\tau(v)$ with $|\tau|<\varepsilon_0$. This proves that the restricted flow on $\mathscr B_s^+(C)$ is positively expansive with reparametrizations. The restricted flow is nonsingular and its phase space is compact. By \cite[Theorem 4.2]{Art14}, $\mathscr B_s^+(C)$ consists of finitely many periodic orbits. Reversing the transverse orientation or the direction of time proves the other cases.
\end{proof}

\begin{proof}[Proof of Proposition~\ref{prop:freeside}]
Let $\varnothing\ne C\subsetneq M$ be compact, invariant, and locally
maximal. Since $M$ is connected, $\partial C\ne\varnothing$.
In a small transverse product rectangle through a boundary point,
write $C\cap R=K^s\times K^u$. At least one factor fails to contain
a neighborhood of the corresponding coordinate. A complementary
interval with an endpoint in that factor, chosen inside the
coordinate interval, therefore gives a stable or unstable gap.
The base of this gap may differ from the original boundary point.
After reversing time or the relevant transverse orientation if
necessary, we obtain a closed positive stable interval $I$ in a
transverse section $R_0$ such that
\[
  I\cap C=\{x\},
\]
where $x$ is the initial endpoint of $I$.

Let $\alpha_{\phi}(x)$ denote the $\alpha$-limit set of the orbit through $x$.
Compactness and invariance give
$\varnothing\ne\alpha_{\phi}(x)\subset C$.
Choose $y\in\alpha_{\phi}(x)$. We claim that $y\in \mathscr B_s^+(C)$.

Suppose otherwise. A point of $\wt C$ lies strictly on the
positive stable side of a lift $\wt y$ of $y$. Choose a fixed
compact transverse arc starting at $\wt y$ and ending on that
orbit in the lifted weak stable plane.
Lemma~\ref{geo:contraction}(3) gives a fixed positive transport time
$T$ and a small terminal product rectangle $R_1$ in whose interior
the transported endpoints project to
\[
  y'=\phi_T(y)=(s_0,u_0),\qquad
  w'=(s_1,u_0)\in C,\qquad s_0<s_1.
\]
Take $r_n\to\infty$ with $\phi_{-r_n}(x)\to y$.
The local flow-box projection to $R_1$ gives numbers $\sigma_n\to0$
such that
\[
  x_n'=\operatorname{pr}_{R_1}(\phi_{T-r_n}(x))=\phi_{T-r_n+\sigma_n}(x)=\phi_{-t_n}(x),
\] where $t_n=r_n-T-\sigma_n\longrightarrow\infty$ and $x_n'=(s_n,u_n)\longrightarrow(s_0,u_0)$.
All these points lie in $C$. For large $n$, local product structure
therefore gives
\[
  z_n=(s_1,u_n)\in C,\qquad s_n<s_1.
\]
Let $J_n\subset R_1^s(x_n')$ be the closed stable interval from
$x_n'$ to $z_n$.

Transport $J_n$ forward along the orbit segment from $x_n'$ to $x$,
using $R_0$ as the terminal section. For large $n$,
Lemma~\ref{geo:contraction}(2), with the local section projections,
gives a map
\[
  H_n(v)=\phi_{\tau_n(v)}(v),\qquad
  \tau_n(x_n')=t_n,\qquad H_n(x_n')=x,
\]
whose image lies in $R_0^s(x)$ and has stable-coordinate width tending
to zero. This is where uniformity over the varying local intervals
$J_n$ is used.
The continued transport preserves the strict order of the lifted
orbit labels. Consequently,
\[
  H_n(z_n)\in I\setminus\{x\}
\]
for large $n$. Invariance also gives $H_n(z_n)\in C$, contradicting
$I\cap C=\{x\}$. Thus $y\in \mathscr B_s^+(C)$.

Lemma~\ref{geo:free-points} now implies that $\mathscr B_s^+(C)$ consists of finitely many periodic orbits. Let $o\subset \mathscr B_s^+(C)$ be any periodic orbit. Let $\wt o$ be one lifted orbit of
$o$, and let $\wt W$ be its weak stable plane.
For a subset $D\subset\wt M$, write $\operatorname{Stab}(D)=\{\gamma\in\operatorname{Deck}(\pi)\mid \gamma(D)=D\}$. Suppose that $\gamma\in\operatorname{Stab}(\wt W)$ does not
preserve $\wt o$. Since $\gamma$ preserves the order of flow
orbits in $\wt W$, either $\gamma\wt o$ or
$\gamma^{-1}\wt o$ lies strictly on the positive side of
$\wt o$. Both orbits belong to $\wt C$, contrary
to freeness. Consequently,
$\operatorname{Stab}(\wt W)=\operatorname{Stab}(\wt o)\cong\mathbb Z$. The last group acts on $\wt o$ by time translations
through integer multiples of the least period of $o$.
Its action on $\wt W$ preserves both the flow orientation and
the transverse orientation. The intrinsic weak leaf $W=\wt W/\operatorname{Stab}(\wt W)$ is therefore an open cylinder, and $o$ is a core circle.
For $\wt x\in\wt o$, the set $\pi(P_s^+(\wt x))$ is one entire component of $W\setminus o$.
Since $P_s^+(\wt x)\cap\pi^{-1}(C)=\varnothing$, we have
$\pi(P_s^+(\wt x))\cap C=\varnothing$. 

The unstable case follows by reversing time.
Undoing the initial relabelling gives the asserted conclusion for
the full preimage $p^{-1}(C)$ in the finite cover. Finiteness follows by taking the union of the four finite families given by Lemma~\ref{geo:free-points}. Therefore, we complete our proof.
\end{proof}

\section{Uniform intersection with the pullback lamination and completion}\label{sec:intersection}

We now show that every orbit of the model flow meets the pullback of an
invariant $su$-lamination. If an orbit were to miss this pullback, the
preceding section would provide a periodic orbit with a free stable or
unstable side. We use that side to construct an essential circle
preserved by an iterate of the orbit equivalence. The collapse map then
sends this circle into a compact invariant subset of one $su$-leaf,
where Corollary~\ref{cor:su-circle} applies.

\begin{lemma}\label{geo:extremal}
Let $S$ be a connected manifold without boundary. Suppose that
$\mathscr A\subset S\times\R$ is a closed codimension-one lamination transverse to the vertical fibers and meeting every fiber. If
$\mathscr A\subset S\times(-\infty,T]$ for some $T\in\R$, then
\[
 \mu(z)=\max\{t\in\R:(z,t)\in \mathscr A\}
\]
is a continuous function, and its graph
\[
 G_\mu=\{(z,\mu(z)):z\in S\}
\]
is a leaf of $\mathscr A$, homeomorphic to $S$. Under a uniform lower
bound, the analogous conclusion holds for the minimum-height graph.
\end{lemma}

\begin{proof}
For each $z\in S$, the set $\{t:(z,t)\in \mathscr A\}$ is nonempty, closed,
and bounded above, so its maximum exists. A graph plaque through any
point of this fiber gives a lower bound for $\mu$ on a neighborhood
of $z$. For $z_n\to z$, this local lower bound and the global upper bound $T$ place $\mu(z_n)$ in a compact interval. Closedness of $\mathscr A$ gives $\limsup_n\mu(z_n)\leq\mu(z)$. On the other hand, the plaque through
$(z,\mu(z))$ is locally the graph of a continuous function $\vartheta$ with
$\vartheta(z)=\mu(z)$. Since $\mu(w)\geq \vartheta(w)$ for nearby $w$, we also have
$\liminf_n\mu(z_n)\geq\mu(z)$. Thus $\mu$ is continuous.

Choose a graph lamination box around $(z,\mu(z))$. By continuity,
the nearby points of $G_\mu$ remain in this box. Let $P$ describe the plaque through $(z,\mu(z))$. A nearby maximal point cannot lie below this plaque, since the plaque itself belongs to $\mathscr A$. If it lay on a
plaque above it, the order of the plaques would give an intersection
with the vertical fiber above $\mu(z)$, which is impossible. Hence
$G_\mu$ locally coincides with a single plaque.

The graph $G_\mu$ is connected and locally a plaque, hence lies in one leaf and is open in that leaf. Since $\mu$ is continuous, its graph is closed in $S\times \R$, and therefore closed in the intrinsic topology of the leaf. It is consequently the whole leaf. Replacing $t$ by $-t$ proves the minimum-height assertion.
\end{proof}

\begin{lemma}\label{geo:missed-set}
Let $\Lambda\ne\varnothing$ be a compact lamination transverse to a
nonsingular flow $\phi$ on a closed manifold $M$. Set
\[
 \mathcal V=\phi_{\R}(\Lambda),\qquad C=M\setminus \mathcal V.
\]
Then $\mathcal V$ is open and invariant, and $C$ is compact and invariant. If
$C\ne\varnothing$, it is locally maximal. If $C=\varnothing$, there
exists $T_0>0$ such that, for every $x\in M$,
\[
 \phi_{(0,T_0)}(x)\cap\Lambda\ne\varnothing,
 \qquad
 \phi_{(-T_0,0)}(x)\cap\Lambda\ne\varnothing.
\]
\end{lemma}

\begin{proof}
The flow saturation of a transverse plaque is open. This proves the
first assertions. If $C\ne\varnothing$, choose a compact neighborhood
$N$ with
\[
 C\subset\operatorname{int}N\subset N\subset M\setminus\Lambda.
\]
An entire orbit contained in $N$ misses $\Lambda$, and hence belongs
to $C$. Conversely, every orbit in $C$ stays in $N$. Thus
$C=\operatorname{Inv}_{\phi}(N)$.

If $\mathcal V=M$, the increasing open sets
\[
 U_n=\{x:\phi_{(-n,n)}(x)\cap\Lambda\ne\varnothing\}
\]
cover $M$. Choose $n_0$ with $U_{n_0}=M$ and apply this equality to
$\phi_{n_0}(x)$ and $\phi_{-n_0}(x)$. The conclusion follows with
$T_0=2n_0$.
\end{proof}

\begin{proposition}\label{prop:invariant-circle}
Let $\Lambda\ne\varnothing$ be a compact lamination transverse to a
topological Anosov flow $\phi$ on a closed three-manifold. Suppose
that a positive self orbit equivalence $\beta$ permutes its leaves
by intrinsic homeomorphisms. If some orbit misses $\Lambda$, then,
after a finite cover, one leaf of $\Lambda$ contains an essential
circle $c$ satisfying $\beta^m(c)=c$ for some $m>0$.
\end{proposition}

\begin{proof}
Let $C=M\setminus\phi_{\R}(\Lambda)$. It is a nonempty proper compact locally maximal invariant set by Lemma~\ref{geo:missed-set}, and $\beta(C)=C$. Pass to the finite orientation cover and a compatible lifted iterate preserving all the orientations.

An orientation-preserving orbit equivalence preserves the weak stable
and weak unstable foliations. Indeed, it takes time-matched
asymptotic orbits to asymptotic orbits, by uniform continuity and the
increasing orbitwise time changes; its inverse gives the reverse
inclusion. After a further iterate it preserves the two weak orders.
It consequently permutes each of the four free-side sets
$\mathscr B_s^\pm(C),\mathscr B_u^\pm(C)$. By Lemma~\ref{geo:free-points} these are
finite unions of periodic orbits, and Proposition~\ref{prop:freeside}
says that at least one is nonempty. We may choose a periodic orbit
$o\subset C$ and one of its free sides $S$ so that
\[
 \beta^m(o)=o,\qquad \beta^m(S)=S
\]
for some $m>0$. Consider a stable side; the unstable case reverses time.

We first construct actual flow coordinates on the intrinsic
half-cylinder $S$. The local return map on its stable side is an
increasing homeomorphism $r$ fixing the origin and satisfying
$0<r(x)<x$ for small $x>0$. It is topologically conjugate near the
origin to $x\mapsto x/2$. In suspension coordinates $(x,\vartheta)$, where
$0\le \vartheta\le1$ and $(x,1)\sim(x/2,0)$, for sufficiently small $\varepsilon_1>0$, the graph $x=\varepsilon_1\,2^\vartheta$ is a
transverse circle bounding a small annular collar of $o$. The flow
crosses it toward the collar. A continuous positive return-time
function changes the speed in these coordinates and preserves this
transversality. Denote the circle by $\Sigma$.

Every orbit in $S$ enters the collar. To see this, join its orbit label
to the periodic axis by a fixed compact stable arc and apply
Lemma~\ref{geo:contraction}. Every nonaxis orbit in the collar leaves
it in finite negative time, since inverse returns multiply its stable
coordinate by $2$ while they remain there. The collar is positively
invariant. Thus every orbit in $S$ crosses $\Sigma$ exactly once.
Flow boxes give a homeomorphism
\[
 \Psi:\Sigma\times\R\longrightarrow S,
 \qquad \Psi(z,t)=\phi_t(z),
\]
with the original time parametrization of $\phi$.

The sets $\phi_t(\Sigma)$ approach $o$ uniformly as $t\to+\infty$.
For a prescribed neighborhood of $o$, choose a smaller trapping collar
inside it. Every point of $\Sigma$ enters its interior. The open sets
of points entering before time $n$ increase to cover compact $\Sigma$;
a common entry time follows, and positive invariance controls all
later times.

Now
\[
 \mathscr A=\Psi^{-1}(\Lambda\cap S)
\]
is a closed one-dimensional lamination transverse to the vertical
fibers, obtained by restricting the graph plaques of $\Lambda$ to
weak stable plaques. Every fiber meets $\mathscr A$, because $S\cap C$ is
empty. The compact disjoint sets $o$ and $\Lambda$ have positive
separation. Uniform approach to $o$ gives
$\mathscr A\subset\Sigma\times(-\infty,T]$ for some $T$.
Lemma~\ref{geo:extremal} gives a circular intersection leaf
\[
 c=\{\Psi(z,\mu(z)):z\in\Sigma\}.
\]
Its local parametrizations lie in source plaques, so it is intrinsically
continuous in one source leaf of $\Lambda$.

On each oriented orbit in $S$, the point of $c$ is its last
intersection with $\Lambda$. This description is independent of
$\Sigma$ and of the flow-time coordinate. Since $\beta^m$ preserves
$S$, $\Lambda$, and the forward order of each orbit, it preserves
these last intersections. Hence $\beta^m(c)=c$.

The circle is a section of $\Sigma\times\R$ and represents the
generator of the weak cylinder. Weak leaves of a topological Anosov
flow are incompressible, so $c$ is essential in $M$. 
First intersections give the same conclusion on an unstable side.
\end{proof}

\begin{proof}[Proof of Theorem~\ref{thm:intersection}]
Let $\Lambda=h^{-1}(\Gamma)$ be the full preimage under the collapse map $h$. By
Proposition~\ref{prop:pullback}, it is a nonempty compact surface
lamination transverse to the model flow, and $\beta$ permutes its
leaves by intrinsic homeomorphisms.

Suppose that some model orbit of $\phi$ misses $\Lambda$.
Proposition~\ref{prop:invariant-circle} gives a finite cover on which one source leaf of $\Lambda$ contains a circle $c$ that is not null-homotopic in the covering manifold and satisfies $\beta^m(c)=c$ for some $m>0$. Retain the same notation for the lifted collapse data and the required iterates, so that $f\circ h=h\circ\beta$.

Proposition~\ref{prop:pullback} also gives intrinsic continuity of $h_c=h|_c:c\longrightarrow L$, where $L$ is the target $su$-leaf containing $h(c)$. Hence $K=h_c(c)$ is a nonempty compact subset of $L$ in its intrinsic topology. With $g=\beta^m|_c$, the collapse relation gives
\[
 f^m\circ h_c=h_c\circ g,
 \qquad f^m(K)=K.
\]
Here $g:c\to c$ is a circle homeomorphism. Since $f$ permutes the
target leaves and $f^m(L)$ meets $L$ in $K$, the two leaves agree.
Corollary~\ref{cor:su-circle} therefore implies that $h_c$ is constant.

This contradicts the essentiality of $c$ and the fact that $h$ is homotopic to the identity. Therefore, every model orbit of $\phi$ must meet $\Lambda$. Applying Lemma~\ref{geo:missed-set} on the original manifold gives the stated uniform positive and negative hitting windows, which descend from the finite cover.
\end{proof}

Once every flow orbit meets the lamination, each complementary flow
segment has two endpoints on the lamination. We next use these
endpoints to describe the complementary regions and extend the
lamination to a foliation.

\begin{proposition}\label{geo:completion}
Let $\Lambda$ be a compact lamination transverse to a nonsingular flow
on a closed three-manifold. Suppose that every orbit meets $\Lambda$.
Each component of $M\setminus\Lambda$ is homeomorphic to
$L\times(0,1)$, where $L$ is an entire boundary leaf and the interval
fibers are flow segments. Adjoining their two endpoint copies gives a
cut-open completion homeomorphic to $L\times[0,1]$.
The lamination extends to a topological foliation of $M$. If its
original leaves are all noncompact, the extension has no compact
leaves.
\end{proposition}

\begin{proof}
For $x\notin\Lambda$, the uniform bound in
Lemma~\ref{geo:missed-set} gives a first preceding and a first following
intersection with $\Lambda$. Their hitting times and endpoint maps
are locally continuous into specified boundary plaques. Indeed, a
compact middle portion of the gap has a neighborhood disjoint from
$\Lambda$. Near either endpoint, the incoming side of its plaque has
no intervening plaque on the central flow fiber. In an ordered graph
chart the same side is empty over a smaller transverse disk. The
continuous hitting graph there remains the first hit.

Let $D^+$ be the set of points of $\Lambda$ having a positive empty
flow side, and define $D^-$ similarly. These sets are unions of whole
leaves. In a lamination box the condition is isolation of the plaque
parameter from the indicated side; holonomy along a finite plaque
chain preserves it. The next-hit map $D^+\to D^-$ and the preceding-hit
map $D^-\to D^+$ are inverse local homeomorphisms in intrinsic
plaques. They pair entire leaves. Each sends a connected leaf into
one leaf; on that whole target leaf the inverse has connected image
in the original leaf, and the two local maps are inverse everywhere.

For a lower boundary leaf $L\subset D^+$ and any point $p\in L$, let $\tau_+(p)>0$ be its
next hitting time. It is intrinsically continuous. The map
\[
 J_L:L\times(0,1)\longrightarrow M\setminus\Lambda,
 \qquad J_L(p,s)=\phi_{s\tau_+(p)}(p),
\]
is injective, because a point in a gap has a unique preceding endpoint
and elapsed time. The local hitting graphs make it a local
homeomorphism. Its image is connected and open. These images are
pairwise disjoint and exhaust the complement, so each is also
relatively closed there. They are precisely its connected components.
Adding the two endpoint copies gives the stated completion, keeping
the boundary sides distinct even when their ambient images lie in the
same leaf.

Fill each complementary component with the surfaces
$J_L(L\times\{s\})$, $0<s<1$. We verify that these surfaces and
the original plaques form foliation charts along $\Lambda$.
Choose actual flow coordinates $(z,t)$ over a transverse disk $D$,
fix $z_0\in D$, and write the original plaques as graphs
$t=\vartheta_r(z)$, indexed by their heights on the central fiber.
Thus $\vartheta_r(z_0)=r$, the parameter $r$ ranges over a relatively
closed set $K$ in an interval, and $(z,r)\mapsto\vartheta_r(z)$ is
jointly continuous. The graphs are strictly ordered by $r$.
If a full flow gap lies in this box between adjacent plaques
$\vartheta_{r_-}<\vartheta_{r_+}$, where $(r_-,r_+)\cap K=\varnothing$,
the inserted graphs have heights
\[
 t=(1-s)\vartheta_{r_-}(z)+s\vartheta_{r_+}(z),\qquad 0\leq s\leq1.
\]
Indeed, the elapsed flow time between the endpoint plaques is
$\vartheta_{r_+}(z)-\vartheta_{r_-}(z)$, and $J_L$ uses this time to parametrize the gap.

For $r=(1-s)r_-+sr_+$, denote the interpolated function by
$\overline\vartheta_r$, and set $\overline\vartheta_r=\vartheta_r$ for $r\in K$.
Continuity within each closed gap follows directly from the
interpolation formula. To check continuity at accumulating gaps,
let $(r_n^-,r_n^+)$ be distinct gaps with
$r_n^-,r_n^+\to r_*\in K$, and choose any $v_n\in[r_n^-,r_n^+]$.
For every compact disk $D_0\subset D$, the ordering gives
\[
 \sup_{z\in D_0}
   |\overline\vartheta_{v_n}(z)-\vartheta_{r_*}(z)|
 \leq
 \max_{r\in\{r_n^-,r_n^+\}}
   \sup_{z\in D_0}|\vartheta_r(z)-\vartheta_{r_*}(z)|
 \longrightarrow0.
\]
The convergence follows from joint continuity of the original
graph family on compact sets. Hence the interpolated graphs
remain jointly continuous across their limiting plaques.

If a full flow gap is not contained in the chosen box, work near
its visible endpoint. At a lower endpoint, choose a fixed
$s_0>0$ sufficiently small. Continuity of the gap-length function
$\tau_+$ ensures, after shrinking $D$, that the corresponding
plaque of $J_L(L\times\{s_0\})$ is a graph inside the box.
Use this graph as an auxiliary boundary. Between it and the
original boundary plaque, the same interpolation formula holds
with parameter $s/s_0$. The upper-endpoint case follows by
reversing time. Including these auxiliary graphs where needed
and shrinking the neighborhood, we obtain an open interval $I$
of height parameters for which
\[
 \Theta:D\times I\longrightarrow M,\qquad
 \Theta(z,r)=(z,\overline\vartheta_r(z)),
\]
is a homeomorphism onto an open neighborhood of the point under
consideration. Indeed, the graph family is jointly continuous
and strictly ordered in $r$. Its horizontal slices are the
plaques of the extended foliation, and its vertical fibers are
flow segments. This proves the required foliation-chart
compatibility along $\Lambda$. Each inserted leaf is homeomorphic
to its boundary leaf $L$, and the original leaves remain leaves
of the extension. Thus the resulting foliation has no compact
leaves whenever all leaves of $\Lambda$ are noncompact.
\end{proof}

\begin{proposition}\label{prop:source-planes}
Let $\Lambda$ be a compact surface lamination transverse to a
nonsingular flow on a closed three-manifold. If every flow orbit meets
$\Lambda$ and $\Lambda$ has no compact leaves, then every lifted leaf
of $\Lambda$ in the universal cover is a properly embedded separating
plane.
\end{proposition}

\begin{proof}
Proposition~\ref{geo:completion} extends $\Lambda$ to a topological
foliation $\mathcal F$ without compact leaves. The flow orients its
transverse direction. It has no Reeb component, since the boundary
of a Reeb component is a compact torus leaf. We may pass to the
orientation cover of the ambient manifold; this leaves the universal
cover unchanged, and the lifted foliation still has no compact leaves.

Novikov's theorem \cite{Nov65} and its extension to topological foliations by Solodov \cite{Sol84} (see also \cite{Cal01} and \cite[Theorem 2.2 and 2.7]{Bow16}) imply that each leaf inclusion induces an injective homomorphism on fundamental groups and that no closed transversal is null-homotopic. The exceptional product foliation by spheres is excluded by the absence of compact leaves.

Let $\pi:\wt M\to M$ be the universal covering projection,
and let $\wt L$ be a leaf of the lifted foliation.
The injectivity just stated implies that
$\pi|_{\wt L}:\wt L\to L$ is the universal covering
of the corresponding leaf $L$. Thus $\wt L$ is simply
connected. It is noncompact, since a compact covering leaf would
have compact image $L$. The classification of surfaces therefore
identifies $\wt L$ with a plane.

To prove properness, suppose that an intrinsic sequence escaping
every compact subset of $\wt L$ converges in $\wt M$.
In a sufficiently small foliation box around its limit, the sequence
must pass through infinitely many distinct plaques. Otherwise,
in a smaller box it would lie in a finite union of compact plaque
disks, contradicting its intrinsic escape. A short transversal in
the box then meets $\wt L$ twice. Join the two intersections
by a path in $\wt L$. The transverse segment and this
leafwise path can be perturbed to a closed transversal.
Its projection to $M$ is null-homotopic because $\wt M$ is
simply connected, which gives a contradiction. Hence the
leaf inclusion is proper, and $\wt L$ is a properly
embedded plane.

Finally, $\wt L$ is two-sided. If its complement were
connected, a short arc crossing it once could be closed by a
path in the complement. This would give a loop with mod-two
intersection number one with $\wt L$. Such an intersection
number is impossible for a null-homotopic loop in a simply
connected manifold. Thus $\wt L$ separates $\wt M$.
The leaves of $\wt\Lambda$ are leaves of this lifted foliation,
so the conclusion applies to them.
\end{proof}

For the accessibility lamination on a manifold whose fundamental
group is not virtually solvable, Lemmas~\ref{lem:no-compact}
and~\ref{lem:no-compact-source} exclude compact leaves of $\Lambda$. Theorem~\ref{thm:intersection} and Proposition~\ref{prop:source-planes} therefore give the properly
embedded source planes needed in the orbit-space argument.

\section{Periodic lozenges and essential circles}\label{sec:lozenge}

The preceding section gives proper lifted planes on the model side of the collapsing map. We use these planes to carry out the lozenge construction, adapting some ideas taken from~\cite[Proposition 4.3]{Fenley99} and~\cite[Section 4 and 5]{FP24}. The resulting circle lies in a source leaf. Its image need not be embedded, so the final obstruction will use Corollary~\ref{cor:su-circle}.

Throughout the geometric part of this section, $\Lambda$ is a nonempty compact lamination transverse to a topological Anosov flow $\phi$ on a closed three-manifold. We assume that every leaf of $\wt\Lambda$ is a properly embedded plane. Write
\[
q:\wt M\longrightarrow\cO=\wt M/\wt\phi.
\]
The space $\wt M$ is homeomorphic to $\R^3$, the orbit space $\cO$ is homeomorphic to $\R^2$, and the lifted flow has a global topological product structure; see~\cite{Fen94,Bar95}. The weak foliations induce foliations $\mathcal G^s,\mathcal G^u$ of $\cO$ by properly embedded lines. We use the same notation for a deck transformation, or a lift of a self orbit equivalence, and its induced action on $\cO$. 

\subsection{The orbit projection of a transverse plane}

\begin{lemma}\label{src:projection}
For every leaf $L$ of $\wt\Lambda$, each lifted flow orbit meets $L$ at most once. The restriction
\[
q|_L:L\longrightarrow U_L=q(L)
\]
is a homeomorphism onto a connected open subset of $\cO$.
\end{lemma}
\begin{proof}
The proper plane $L$ separates $\wt M$. Its transverse flow orientation specifies a positive side, and every positive-time crossing goes from the negative side to the positive side. Two crossings would require an intervening crossing in the opposite direction. Thus an orbit meets $L$ at most once. Transverse flow boxes make $q|_L$ a local homeomorphism. It is injective and therefore a homeomorphism onto its open image, which is connected.
\end{proof}

Compactness of the lamination gives a uniform form of transversality. For any $\eps>0$, there exists $\delta>0$ such that
\begin{equation}\label{src:uniform-crossing}
x\in L\in\wt\Lambda,\quad d(x,y)<\delta
\quad\Longrightarrow\quad
\wt\phi_{[-\eps,\eps]}(y)\cap L\ne\varnothing.
\end{equation}
Indeed, cover $\Lambda$ by finitely many smaller flow boxes with closures in transverse lamination boxes. Each plaque in a larger box is a graph over its transverse disk. Nearby flow fibers cross the specified graph, and compactness of the smaller boxes gives uniform spatial and time margins. Lift these boxes to $\wt M$. In particular, the intersection in~\eqref{src:uniform-crossing} belongs to the specified leaf $L$.

We give the boundary argument in the topological setting needed here. It is the source version of the argument in~\cite{FP24}.

\begin{lemma}\label{src:boundary}
Let $L\in\wt\Lambda$ and $U=q(L)$. If $o\in\partial U$, exactly one of
\[
\mathcal G^s(o)\subset\partial U,
\qquad
\mathcal G^u(o)\subset\partial U
\]
holds. The set $U$ lies on one side of the corresponding line.
\end{lemma}
\begin{proof}
Write $\wt o=q^{-1}(\{o\})$ for the corresponding lifted orbit.
Take a small transverse product rectangle $D$ through $\wt o$ and identify it with its image in $\cO$. For $v\in D\cap U$, let $\tau(v)$ be the time from its representative in $D$ to its intersection with $L$. Uniqueness and transversality make $\tau$ continuous on its domain. If $v_n\in D\cap U$ tends to $o$, then
\[
|\tau(v_n)|\longrightarrow\infty.
\]
Otherwise the intersection points have a convergent subsequence on the lifted orbit $\wt o$. Since $L$ is properly embedded and hence closed in $\wt M$, the orbit $\wt o$ would meet $L$.

Suppose, after passing to a subsequence, that $\tau(v_n)\to+\infty$. Uniform local stable contraction implies that a fixed small stable plaque through $v_n$ is contained in $U$ for all large $n$. Here and below a stable plaque is the specified plaque in $D$, not another return of the same global weak leaf. To see the assertion, continue the plaque by positive stable Poincar\'e holonomy to the intersection point of $v_n$ with $L$. The terminal plaques become arbitrarily small by Lemma~\ref{geo:contraction}; equivalently, one may apply~\cite[Lemmas 5--7]{Iak25}. Equation~\eqref{src:uniform-crossing} then makes all their orbits cross the same local plaque of $L$ within a uniformly bounded additional time.

We also need the corresponding escape statement. Set
\[
x_n=\wt\phi_{\tau(v_n)}(v_n)\in L.
\]
These points escape every compact subset of $\wt M$, since a convergent subsequence would have its limit in $L$ over the orbit $o$. If $w_n$ belongs to the specified stable plaque of $v_n$, its transported point is arbitrarily close to $x_n$, and its intersection with $L$ differs from that point by a uniformly bounded flow time. Its intersections therefore escape compact subsets as well. The holonomy has been continued in positive time from $D$; its terminal times tend to $+\infty$, since bounded nonnegative times would keep the terminal points in a compact set. The bounded additional crossing times do not change this conclusion. Thus
\[
w_n\in U,\qquad \tau(w_n)\longrightarrow+\infty.
\]
This argument uses local Poincar\'e contraction and makes no assertion about a uniform global time change along a weak leaf.

For any $w$ on the local stable plaque of $o$, choose $w_n$ on the stable plaque of $v_n$ with $w_n\to w$. The fact $w_n\in U$ shows that $w\in\overline U$. If $w\in U$, continuity of its finite hitting time would contradict that equation. Hence the entire local stable plaque of $o$ lies in $\partial U$.

Any prescribed compact segment of $\mathcal G^s(o)$ is covered by a finite chain of transverse product rectangles. The positive-diverging hitting sequence in one rectangle supplies such a sequence in the next; changing representatives on their overlap changes time by a bounded amount. Repeating the argument along this chain gives
\[
\mathcal G^s(o)\subset\partial U.
\]
Since $U$ is connected and this line separates $\cO$, the set $U$ lies on one side of it. The local stable plaques through $v_n$ cross $\mathcal G^u(o)$ and lie in $U$. Thus $\mathcal G^u(o)$ cannot also be contained in $\partial U$.

If the negative-diverging alternative $\tau(v_n)\longrightarrow -\infty$ occurs, the same proof uses unstable plaques and gives the unstable boundary line. The preceding exclusivity also shows that the two alternatives cannot arise from different subsequences.
\end{proof}

\subsection{A canonical separating leaf}

Recall that $q:\wt M\to\cO$ is the quotient map by the
orbits of the lifted flow $\wt\phi$, and that
$\mathcal G^s$ and $\mathcal G^u$ are the foliations of $\cO$
induced by the lifted weak stable and weak unstable foliations.
A half-leaf based at $o\in\cO$ is a connected component of
$\mathcal G^\sigma(o)\setminus\{o\}$, where $\sigma=s$ or $u$.

We use the conventions of~\cite{Fen94} and the orbit-space
description in~\cite[Section~3]{Fenley98}.
An embedding of a subset of $\mathbb R^2$ into $\cO$ is
\emph{foliated} if it sends horizontal and vertical segments
into stable and unstable leaves, respectively.
It is \emph{proper} if the inverse image of every compact
subset of $\cO$ is compact.

Disjoint leaves $\ell^s\in\mathcal G^s$ and $\ell^u\in\mathcal G^u$
form a \emph{perfect fit} if there is a proper foliated embedding
\[
  \iota:[0,1]^2\setminus\{(1,1)\}\longrightarrow\cO
\]
such that
\[
  \iota([0,1)\times\{1\})\subset \ell^s,
  \qquad
  \iota(\{1\}\times[0,1))\subset \ell^u.
\]
Thus the two leaves bound a foliated rectangle with one
missing corner. Properness means that approaching this
missing corner in the model forces the image to leave every
compact subset of $\cO$.
We also say that a stable half-leaf and an unstable half-leaf
form a perfect fit when they contain the two boundary rays
in such a model.

A \emph{lozenge} is an open region
$\mathcal L=\xi((0,1)^2)\subset\cO$, where
\[
  \xi:D\longrightarrow\cO,
  \qquad
  D=[0,1]^2\setminus\{(0,1),(1,0)\},
\]
is a proper foliated embedding. Its corners are
$o_1=\xi(0,0)$ and $o_2=\xi(1,1)$. Choose the orientations of $\mathcal G^s$ and $\mathcal G^u$ so that the two sides of $\mathcal L$ based at $o_1$ are
the positive half-leaves. In the coordinates above, the four
sides are then
\[
  \begin{aligned}
    \mathcal G^s_+(o_1)&=\xi((0,1)\times\{0\}),&
    \mathcal G^u_+(o_1)&=\xi(\{0\}\times(0,1)),\\
    \mathcal G^s_-(o_2)&=\xi((0,1)\times\{1\}),&
    \mathcal G^u_-(o_2)&=\xi(\{1\}\times(0,1)).
  \end{aligned}
\]
Here $\mathcal G^\sigma_\pm(o)$ denotes the positive or
negative component of
$\mathcal G^\sigma(o)\setminus\{o\}$, for $\sigma=s,u$.
In particular,
\[
  \partial\mathcal L
  =\{o_1,o_2\}
   \cup\mathcal G^s_+(o_1)\cup\mathcal G^u_+(o_1)
   \cup\mathcal G^s_-(o_2)\cup\mathcal G^u_-(o_2).
\]
The pairs
\[
  \bigl(\mathcal G^s_+(o_1),\mathcal G^u_-(o_2)\bigr),
  \qquad
  \bigl(\mathcal G^s_-(o_2),\mathcal G^u_+(o_1)\bigr)
\]
form the two perfect fits.
The missing vertices of $D$ do not correspond to points of $\cO$.

Every stable leaf meeting $\mathcal L$ intersects both
unstable sides, and every unstable leaf meeting $\mathcal L$
intersects both stable sides, in each case exactly once.
These intersections define homeomorphisms
\[
  \mathcal H^u:\mathcal G^s_+(o_1)\longrightarrow\mathcal G^s_-(o_2),
  \qquad
  \mathcal H^s:\mathcal G^u_+(o_1)\longrightarrow\mathcal G^u_-(o_2),
\]
characterized by
\[
  \{\mathcal H^u(p)\}=\mathcal G^u(p)\cap\mathcal G^s_-(o_2),
  \qquad
  \{\mathcal H^s(r)\}=\mathcal G^s(r)\cap\mathcal G^u_-(o_2).
\]
When each half-leaf is ordered away from its corner,
both maps reverse order.

A lozenge is \emph{periodic} if
$\gamma(\mathcal L)=\mathcal L$ for some nontrivial deck
transformation $\gamma$.
After replacing $\gamma$ by a positive power, we may assume
that it fixes both corners and preserves each of their stable
and unstable half-leaves.
Write
\[
  \wt o_i=q^{-1}(\{o_i\})
\]
for the corresponding lifted flow orbits.
We also write $W_i^\sigma=q^{-1}(\mathcal G^\sigma(o_i))$ for the
lifted weak plane through $\wt o_i$, where $\sigma=s,u$.
After replacing $\gamma$ by its inverse if necessary, there
are $T_1,T_2>0$ such that
\[
  \gamma(x_1)=\wt\phi_{T_1}(x_1)
  \quad(x_1\in\wt o_1),
  \qquad
  \gamma(x_2)=\wt\phi_{-T_2}(x_2)
  \quad(x_2\in\wt o_2).
\]
The opposite signs are a property of periodic lozenges;
see~\cite{Fen94} and also~\cite[Proposition~2.4]{FP24}.
The two positive numbers $T_1$ and $T_2$ need not be equal.

For each $i\in\{1,2\}$, $\sigma\in\{s,u\}$, and
$\varepsilon\in\{+,-\}$, the restriction
\[
  \gamma|_{\mathcal G^\sigma_\varepsilon(o_i)}
  :\mathcal G^\sigma_\varepsilon(o_i)
   \longrightarrow\mathcal G^\sigma_\varepsilon(o_i)
\]
is an increasing homeomorphism without fixed points.
Its positive and negative iterates approach the two different
ends of this half-leaf, one of which is $o_i$.
To state this explicitly, order
$\mathcal G^\sigma_\varepsilon(o_i)$ away from $o_i$,
and denote this order by $\prec$.
For every $p,x_-,x_+\in\mathcal G^\sigma_\varepsilon(o_i)$
with $x_-\prec x_+$, there exist integers $m,n$ such that
\[
  \gamma^m(p)\prec x_-\prec x_+\prec\gamma^n(p).
\]
This is the cofinality of the iterates toward both ends
that we use below.
Indeed, a monotone sequence of iterates with a limit in
the interior of the half-leaf would have a fixed point
as its limit.

\begin{proposition}\label{src:frontier}
Suppose that $\wt o_1$ meets $\wt\Lambda$. There is a $\gamma$-invariant leaf $L$ of $\wt\Lambda$ which misses both lifted corner orbits and meets the lift of an orbit in the interior of $\mathcal L$. Moreover, if a lift $\wt b$ of a positive self orbit equivalence $b$ permutes the leaves of $\wt\Lambda$ and preserves $\mathcal L$, then a positive power of $\wt b$ preserves $L$ and commutes with $\gamma$.
\end{proposition}
\begin{proof}
Choose a leaf $F$ meeting $\wt o_1$ and put $F_n=\gamma^n F$. These leaves are distinct. Otherwise a nontrivial power of $\gamma$ would fix $F$ and give two intersections of $F$ with $\wt o_1$, contrary to Lemma~\ref{src:projection}. Their intersections with $\wt o_1$ strictly order the planes. Let $V_n$ be the open slab between $F_n$ and $F_{n+1}$, and define
\[
\mathcal R(F)=\bigcup_{n\in\Z}(V_n\cup F_n).
\]
The two adjacent slabs fill both local sides of $F_n$, so $\mathcal R(F)$ is open. It is connected and $\gamma$-invariant. The deck transformation acts by a nonzero constant time translation on $\wt o_1$, and hence the intersections $F_n\cap \wt o_1$ are cofinal in both directions. Therefore $\wt o_1\subset \mathcal R(F)$.

The region is independent of the chosen leaf meeting $\wt o_1$. Indeed, if $F'$ is another such leaf, then for some integer $j$ it equals $F_j$ or lies in $V_j$. In the latter case separation and disjointness give
\[
F_{j+n}<\gamma^n F'<F_{j+n+1},\qquad n\in\Z,
\]
where the inequalities refer to the nested sides determined by the flow orientation. Every slab bounded by finitely many planes of either sequence is contained between two planes of the other sequence. Taking their increasing unions yields
\[
\mathcal R(F')=\mathcal R(F).
\]

The orbit $\wt o_2$ is disjoint from all $F_n$. Otherwise $F_{n+1}$ would lie on the positive side of $F_n$ when read along $\wt o_1$, and on its negative side when read along $\wt o_2$. This is impossible for disjoint separating planes. The orbit $\wt o_2$ cannot meet $V_n$ either. Indeed, if $x\in \wt o_2\cap V_n$, then $\gamma x\in V_{n+1}$ precedes $x$ in flow time, and the positive orbit segment from $\gamma x$ to $x$ crosses $F_{n+1}$ in the wrong direction. Consequently
\[
\wt o_2\cap \mathcal R(F)=\varnothing.
\]

We now describe the frontier of $\mathcal R(F)$. Every frontier point is accumulated by plaques of $F_n$ with unbounded indices. Otherwise a small ball missing these planes would have constant slab membership. Properness excludes accumulation of a fixed finite family. In a transverse lamination box, each vertical flow fiber meets each $F_n$ at most once. The approaching plaques are therefore ordered graphs converging monotonically to the frontier plaque. Between two approaching plaques, separation forces every intermediate $F_n$ to cross the intervening vertical segments; these segments belong to the slabs in the region $\mathcal R(F)$. After shrinking the box, the entire local side toward $\mathcal R(F)$ of the limiting plaque is filled by $\mathcal R(F)$. It is clear that each frontier point is contained in $\wt\Lambda$ by the compactness.

Continuation through a finite chain of lamination boxes transports the same limiting plaques along any compact path in their limiting leaf. Thus the frontier is a union of whole leaves of $\wt\Lambda$. For a frontier leaf $L'$, denote the component of $\wt M\setminus L'$ containing $\mathcal R(F)$ by $I_{L'}$, and the other component by $O_{L'}$. The local $I_{L'}$-side is filled by $\mathcal R(F)$. Any orbit crossing $L'$ therefore enters $\mathcal R(F)$ on one of its local sides. The fact $\wt o_2\cap \mathcal R(F)=\varnothing$ implies that $\wt o_2$ misses the frontier, and hence
\[
\wt o_2\cap\overline{\mathcal R(F)}=\varnothing.
\]

Take a path from a point of $\wt o_2$ to $\mathcal R(F)$ and let its first contact with $\overline{\mathcal R(F)}$ belong to the leaf $L\in \wt\Lambda$. Immediately before contact it lies in $O_L$, since the other local side is filled by $\mathcal R(F)$. Its earlier portion avoids $L$, and $\wt o_2$ is connected and disjoint from $L$. It follows that $\wt o_2\subset O_L$.

This frontier leaf is unique. Suppose that $L'$ is another frontier leaf with $\wt o_2\subset O_{L'}$, then
\[
L'\subset I_L,\qquad L\subset I_{L'},
\]
because both lie in $\overline{\mathcal R(F)}$ and are disjoint. The connected set $O_L$ avoids $L'$ and has boundary on $L\subset I_{L'}$, so $O_L\subset I_{L'}$. It follows that $O_L\cap O_{L'}=\varnothing$. However, both $O_L$ and $O_{L'}$ contain $\wt o_2$, which gives uniqueness of the frontier leaf $L$. Since $\gamma$ preserves $\mathcal R(F)$ and $\wt o_2$, it preserves this unique leaf $L$.

The leaf $L$ is disjoint from both $\wt o_1\subset \mathcal R(F)$ and $\wt o_2$. A path in $\overline{\mathcal L}$ joining the two corners $o_1$ and $o_2$, with all other points in its interior, lifts to a path in $\wt M$ by the orbit product structure. Its endpoints lie on opposite sides of $L$, so it meets $L$ above the interior of $\mathcal L$.

Now assume that $\wt b$ is a lift of a positive self orbit equivalence preserving $\wt \Lambda$ and $\mathcal L$. Replace it by a square if necessary so that its induced action on orbit space fixes both corners and the relevant half-leaves. The deck stabilizer of $\wt o_1$ is infinite cyclic. Indeed, its elements act by constant flow-time translations, and proper discontinuity makes their translation times a discrete subgroup of $\R$. Freeness makes the translation homomorphism injective. Let $\kappa$ be its positive generator.

Conjugation by $\wt b$ preserves this stabilizer. Since $\wt b$ preserves the orientation of $\wt o_1$, the conjugate $\wt b\kappa \wt b^{-1}$ is its positive generator, hence equals $\kappa$. The element $\gamma$ is a power of $\kappa$, so $\wt b\gamma=\gamma \wt b$. The leaf $\wt b(F)$ meets $\wt o_1$. Since $\wt b(\wt o_1)=\wt o_1$ and the region $\mathcal R(F)$ is independent of the chosen leaf meeting $\wt o_1$, we have
\[
\wt b(\mathcal R(F))=\mathcal R(\wt b(F))=\mathcal R(F).
\]
The uniqueness of the frontier leaf separating this region from $\wt o_2$ implies $\wt b(L)=L$. Hence, we finish the proof.
\end{proof}

\subsection{The source circle}

\begin{proposition}\label{src:circle}
Under the hypotheses of Proposition~\ref{src:frontier}, there is a properly embedded line $\ell \subset L$, contained in a weak stable or weak unstable plane through $\wt o_1$, such that
\[
\gamma \ell=\ell,\qquad \ell_0=\ell/\langle\gamma\rangle\cong S^1.
\]
The projection of $\ell_0$ to $M$ represents the nontrivial conjugacy class of $\gamma$, up to reversal. Under the additional hypothesis on $\wt b$ in Proposition~\ref{src:frontier}, a positive power of $\wt b$ preserves $\ell$ and induces a homeomorphism of $\ell_0$.
\end{proposition}
\begin{proof}
Let $U=q(L)$. By Lemma~\ref{src:projection} and Proposition~\ref{src:frontier},
\[
\gamma U=U,\qquad U\cap\mathcal L\ne\varnothing,
\qquad o_1,o_2\notin U,
\]
and $U$ is open and connected. We first show that $\mathcal L\subset U$, as in~\cite{FP24}. If $z\in\mathcal L\cap\partial U$, Lemma~\ref{src:boundary} puts a full stable or unstable line through $z$ in $\partial U$. Suppose it is the unstable line $\mathcal G^u(z)$. Every $\mathcal G^u$-leaf crossing $\mathcal L$ intersects the half $\mathcal G^s$-leaf of $o_1$, which parametrizes all crossing $\mathcal G^u$ leaves. The action of $\gamma$ on the given half leaf of $\mathcal G^s(o_1)$ is fixed-point free. It follows that the iterates of $\gamma$ on this half line are cofinal toward both ends. Then, for any $v\in U\cap\mathcal L$, the points $\gamma^nv\in U$ occur on both sides of $\mathcal G^u(z)$. This contradicts connectedness of $U$. The stable alternative is identical. Hence $\mathcal L\cap\partial U=\varnothing$, and connectedness of $\mathcal L$ implies the conclusion $\mathcal L\subset U$.

Both corners $o_1$ and $o_2$ belong to $\partial U$. Suppose that the boundary alternative at the first corner is
\[
\mathcal G^s(o_1)\subset\partial U.
\]
Let $\mathcal G^u_+(o_1)$ be the unstable half-leaf of $o_1$ lying on the boundary of $\mathcal L$. Lemma~\ref{src:boundary} implies $\mathcal G^u_+(o_1)\subset\overline U$. If $z\in \mathcal G^u_+(o_1)\setminus U$, then $z\in\partial U$. Its stable line $\mathcal G^s(z)$ meets the interior of $\mathcal L$ and cannot be a boundary line. Lemma~\ref{src:boundary} implies that $\mathcal G^u(z)=\mathcal G^u(o_1)\subset \partial U$, contrary to the uniqueness of the boundary alternative at $o_1$. Thus $\mathcal G^u_+(o_1)\subset U$. The other half of $\mathcal G^u(o_1)$ lies across the separating line $\mathcal G^s(o_1)$, so
\[
U\cap\mathcal G^u(o_1)=\mathcal G^u_+(o_1).
\]
It follows that
\[
\ell=(q|_L)^{-1}(\mathcal G^u_+(o_1))
  =L\cap W_1^u,
\]
which is an embedded line. It is closed in $L$ with respect to the intrinsic topology, because $\mathcal G^u_+(o_1)$ is the intersection of $U$ with a full unstable line. It is consequently proper in $L$ and in $\wt M$. 

Both $L$ and $\mathcal G^u_+(o_1)$ are $\gamma$-invariant, which gives $\gamma \ell=\ell$. A fixed-point-free increasing homeomorphism of a line is conjugate to a translation, so the quotient $\ell_0=\ell/\langle\gamma\rangle$ is a circle. A fundamental segment joins a point $x$ to $\gamma x$, and its projection represents $\gamma\ne\id$.

Up to a double power, the map $\wt b$ given in Proposition~\ref{src:frontier} preserves $L$, the corner $o_1$ in orbit space, and the lozenge $\mathcal L$. It therefore preserves the specified half-leaf $\mathcal G^u_+(o_1)$ and thus preserves $\ell$. Since $\wt b\gamma=\gamma \wt b$, it descends to a homeomorphism of $\ell_0$. The other alternative $\mathcal G^u(o_1)\subset\partial U$ follows symmetrically.
\end{proof}

The circle $\ell_0$ is an abstract quotient; its projection to $M$ may have self-intersections. Its map into the intrinsic source leaf is nevertheless continuous, as is seen in the local plaque charts. This distinction is needed when applying the collapsing map.

\begin{proposition}\label{prop:invariant-lozenge}
Let $f$ be a three-dimensional $C^1$ collapsed Anosov flow with data $(\phi,h,\beta)$. Let $\Gamma$ be a nonempty compact $f$-invariant $su$-lamination. Suppose that every lifted leaf of $h^{-1}(\Gamma)$ is a properly embedded plane. Then, for any periodic lozenge $\mathcal L$ having a corner $o_1$ whose lifted orbit $\wt o_1$ meets the lifted lamination of $h^{-1}(\Gamma)$, any lift of a positive iterate of $\beta$ cannot preserve $\mathcal L$.
\end{proposition}
\begin{proof}
Suppose that such a lift exists. Proposition~\ref{prop:pullback} makes $\Lambda=h^{-1}(\Gamma)$ a transverse lamination, with its leaves permuted intrinsically by $\beta$. Apply Proposition~\ref{src:circle}, taking a further positive power of the lift as required. We obtain a source circle $\ell_0$, a circle homeomorphism $g:\ell_0\to \ell_0$, and an integer $m>0$. If $\pi:\wt M\to M$ is the covering projection, the map
\[
\alpha:\ell_0\longrightarrow M,\qquad \alpha([x])=h(\pi(x)),
\]
satisfies
\[
f^m\alpha=\alpha g.
\]
The projection of $\ell_0$ lies intrinsically in one source leaf of $\Lambda$, and Proposition~\ref{prop:pullback} implies that $\alpha$ is intrinsically continuous into one target $su$-leaf $F\in\Gamma$. The set $K=\alpha(\ell_0)$ is intrinsically compact, nonempty, and $f^m$-invariant. Since $f$ preserves the lamination $\Gamma$, this also gives $f^m(F)=F$.

Corollary~\ref{cor:su-circle}, applied to the relation $f^m\alpha=\alpha g$, makes $\alpha$ constant. But the loop $\ell_0$ represents $\gamma\ne\id$, and $h$ is homotopic to the identity. Its image cannot be null-homotopic. This is a contradiction.
\end{proof}

For strong collapsed Anosov flows, the weak-image surfaces give the contradiction without requiring an invariant lozenge for the self orbit equivalence.

\begin{proposition}\label{prop:strong-corner}
Let $f$ be a three-dimensional $C^1$ strong collapsed Anosov flow with data $(\phi,h,\beta)$. Let $\Gamma$ be a nonempty compact $f$-invariant $su$-lamination. Suppose that every lifted leaf of $h^{-1}(\Gamma)$ is a properly embedded plane. Then, no lifted corner orbit of a periodic model lozenge meets the lifted lamination of $h^{-1}(\Gamma)$.
\end{proposition}
\begin{proof}
Suppose otherwise and apply Proposition~\ref{src:circle}, without its additional equivariance conclusion. Let $\wt h$ be the deck-equivariant lift of $h$ supplied by its identity homotopy. The line $\ell$ lies in one lifted weak stable or weak unstable plane $W$. Strong collapse makes $\wt h(W)$ a properly embedded $C^1$ surface tangent to $E^{cs}$ or $E^{cu}$; see~\cite[Definition 2.10 and Lemma 8.2]{BFP23}. By Proposition~\ref{prop:pullback}, the image $\wt h(\ell)$ also lies intrinsically in the lifted target $su$-leaf $F\in \wt \Gamma$ containing $\wt h(L)$, where $L\in \wt \Lambda$ is the leaf given by Proposition~\ref{src:frontier}.

Locally the two specified surfaces $F$ and $\wt h(W)$ intersect in a curve tangent to $E^s$ or $E^u$. Proper embedding of $\wt h(W)$ gives intrinsic continuity into that surface, while the pullback plaque property gives intrinsic continuity into the $su$-leaf. Following their local intersections along the connected line $\ell$ puts $\wt h(\ell)$ in one strong leaf $\ell^{\mathrm{str}}$. Since $\gamma \ell=\ell$ and $\wt h$ is deck-equivariant, the nonempty set $\wt h(\ell)$ is $\gamma$-invariant. Uniqueness of the strong leaf through each point implies $\gamma \ell^{\mathrm{str}}=\ell^{\mathrm{str}}$, contradicting Lemma~\ref{lem:strong-lines}.
\end{proof}

\section{Accessibility and suspension models}\label{sec:main}

We first apply the source construction when strong collapse is available. This will also settle the two exceptional alternatives in the general case.

\begin{proposition}\label{prop:strong-accessibility}
Let $f$ be a $C^1$ strong collapsed Anosov flow on a closed connected three-manifold with non-virtually-solvable fundamental group. Then $f$ is accessible.
\end{proposition}
\begin{proof}
If $f$ were not accessible, $\Gamma(f)$ would be a nonempty invariant $su$-lamination without compact leaves, by Lemma~\ref{lem:no-compact}. Its pullback $\Lambda=h^{-1}(\Gamma)$ has no compact leaves by Lemma~\ref{lem:no-compact-source}. Theorem~\ref{thm:intersection} gives uniform intersection with the pullback lamination, and Proposition~\ref{prop:source-planes} gives proper plane lifts. The nonsuspension model has a periodic lozenge~\cite[Corollary E]{Fenley98}\cite[Theorem 3.3]{Fenley95}. Its corner orbits meet $\wt\Lambda$, contradicting Proposition~\ref{prop:strong-corner}.
\end{proof}

\subsection{The general case}
We recall the following classification result. 

\begin{proposition}\cite[Proposition 2.5]{FP24}\label{FP-3cases}
    Let $\beta$ be a self orbit equivalence of a topological Anosov flow on a three-manifold with non-virtually-solvable fundamental group. Then, one of the following alternatives occurs:
    \begin{enumerate}
\item A lift of a positive iterate of $\beta$ preserves a periodic lozenge.
\item The ambient manifold is hyperbolic.
\item The ambient manifold is Seifert fibered and $\beta$ induces a pseudo-Anosov mapping class on its base.
\end{enumerate}
\end{proposition}

The hyperbolic alternative in that classification has additional information about the model flow and the orbit equivalence, which is unnecessary here. This classification has no transitivity or nonwandering assumption.

We also use the strong-collapse conclusions of Fenley and Potrie~\cite{FPP24}. Partially hyperbolic diffeomorphisms of closed hyperbolic three-manifolds admit strong-collapse data~\cite[Theorem B and Remark 1.2]{FPP24}; the same holds on Seifert manifolds with hyperbolic orientable base and pseudo-Anosov base action~\cite[Theorem C and Remark 1.2]{FPP24}. These conclusions include the additional weak-surface properties of strong collapse and impose no recurrence assumption. The data they supply may differ from the original $(\phi,h,\beta)$.

\begin{proof}[Proof of Theorem~\ref{thm:main}]
Suppose that $f$ is not accessible. Passing to a finite lift of a finite iterate preserves the contradiction argument, so we still use the notation $f$ and its data $(\phi,h,\beta)$. The lamination $\Gamma=\Gamma(f)$ is nonempty and has no compact leaves by Lemma~\ref{lem:no-compact}. By Lemma~\ref{lem:no-compact-source}, Theorem~\ref{thm:intersection}, and Proposition~\ref{prop:source-planes}, its pullback
\[
\Lambda=h^{-1}(\Gamma)
\]
meets every model orbit of $\phi$ and has proper plane lifts on the universal cover $\wt M$.

In the first case of Proposition~\ref{FP-3cases}, the corner orbits of the invariant periodic lozenge meet $\wt\Lambda$. Proposition~\ref{prop:invariant-lozenge} gives a contradiction.

In the second alternative of Proposition~\ref{FP-3cases}, the strong-collapse classification applies to $f$, and Proposition~\ref{prop:strong-accessibility} gives a contradiction. In the Seifert alternative of Proposition~\ref{FP-3cases}, pass if necessary to a finite cover with orientable hyperbolic surface base and lift a positive iterate. The base mapping class remains pseudo-Anosov. Since $h\simeq\id$ and $fh=h\beta$, the maps $f$ and $\beta$ induce the same outer automorphism of $\pi_1(M)$, hence the same base mapping class. The Seifert strong-collapse conclusion applies to this lifted iterate. Its fundamental group is still not virtually solvable, so Proposition~\ref{prop:strong-accessibility} again gives accessibility. Projecting strong paths contradicts the nonaccessibility of the original map.

These alternatives exhaust the possibilities. Hence, we finish the proof.
\end{proof}

The use of strong data in the last two alternatives does not identify general and strong collapsed Anosov flows. The first alternative has been treated using only the original general collapse.

\subsection{Suspension models}
For the corollary, preservation of every model orbit must be upgraded to a continuous time displacement. We record the elementary step that provides this upgrade.

\begin{lemma}\label{lem:continuous-time}
Let $\phi$ be a nonsingular flow on a compact manifold. Suppose that $J_s$, $0\le s\le1$, is an isotopy of homeomorphisms with $J_0=\id$ and that every $J_s$ preserves each $\phi$-orbit. There is a continuous function $\tau:M\times[0,1]\to\R$ satisfying
\[
J_s(x)=\phi_{\tau(x,s)}(x),\qquad \tau(x,0)=0.
\]
\end{lemma}
\begin{proof}
A finite flow-box cover gives a positive lower bound $p_0$ for the periods of nontrivial periodic orbits. Choose smaller flow boxes whose time heights are less than $p_0/4$. Uniform continuity and compactness give a subdivision
\[
0=s_0<s_1<\cdots<s_r=1
\]
such that, for each $x$ and $i$, the path $J_s(x)$, $s_i\le s\le s_{i+1}$, lies in one such box. The same box works for all points near $x$.

An individual orbit meets a flow box in at most countably many plaques. The components of its open time preimage are disjoint open intervals, each containing a rational number. The transverse coordinate of $s\mapsto J_s(x)$ is continuous and takes values in this countable set, so it is constant. Subtracting flow coordinates gives a local continuous increment $\tau_i(x,s)$ with
\[
J_s(x)=\phi_{\tau_i(x,s)}J_{s_i}(x),
\qquad \tau_i(x,s_i)=0.
\]
Two increments on an overlap differ by a period and have difference smaller than $p_0$, hence agree. Their finite sums define the required jointly continuous function $\tau$.
\end{proof}

\begin{proposition}\label{prop:suspension}
Let $f$ be a $C^1$ collapsed Anosov flow whose model flow $\phi$ is orbit equivalent to a suspension Anosov flow. Then, some positive iterate of $f$ is a discretized Anosov flow.
\end{proposition}
\begin{proof}
The model flow $\phi$ is orbit equivalent to the suspension of a hyperbolic toral automorphism~\cite[Proposition 2.7]{BF13}. Conjugate $\beta$ by this orbit equivalence. A positive iterate of a homeomorphism of the resulting mapping torus is homotopic to the identity~\cite[Proposition 1.3]{HP15}. The suspension case of~\cite[Theorem 1.1]{BG19}, together with Remark 1.2 there, gives an isotopy from the identity to that iterate through homeomorphisms preserving each suspension orbit. The transitivity assumption of that result holds for this toral suspension.

Conjugating the isotopy back and applying Lemma~\ref{lem:continuous-time}, we obtain
\[
\beta^m(x)=\phi_{\tau(x)}(x)
\]
with $\tau$ continuous. Thus $\beta^m$ is a trivial self orbit equivalence of $\phi$. The same collapse satisfies $f^mh=h\beta^m$. The proof of~\cite[Proposition 5.27]{BFP23} shows that a general collapsed Anosov flow with trivial associated orbit equivalence is a discretized Anosov flow. Therefore, the diffeomorphism $f^m$ is a discretized Anosov flow.
\end{proof}

\begin{proof}[Proof of Corollary~\ref{cor:nondaf}]
A suspension model would make a positive iterate of $f$ a discretized Anosov flow by Proposition~\ref{prop:suspension}. The hypothesis therefore forces a nonsuspension model. By~\cite[Theorem 5.7]{BF13}, the fundamental group is not virtually solvable, and Theorem~\ref{thm:main} applies.
\end{proof}

\newcommand{\setbiblabel}[2]{%
  \expandafter\def\csname biblabel@#1\endcsname{#2}}
\let\originalbibitem\bibitem
\renewcommand{\bibitem}[2][]{%
  \ifcsname biblabel@#2\endcsname
    \originalbibitem[\csname biblabel@#2\endcsname]{#2}%
  \else
    \originalbibitem[#1]{#2}%
  \fi}

\setbiblabel{BRHTU08}{BHH{\etalchar{+}}08}
\setbiblabel{CHHU18}{CHHU18}
\setbiblabel{HHTU11}{HHTU11}
\setbiblabel{HHTU12}{HHTU12}
\setbiblabel{HHU08}{HHU08a}
\setbiblabel{HHU08b}{HHU08b}
\setbiblabel{HHU11}{HHU11}
\setbiblabel{HHU16}{HHU16}

\bibliographystyle{alpha}
\bibliography{ref}

\end{document}